\documentclass[oneside, reqno] {amsart}

\usepackage{xypic}
\input xy
\xyoption{all}
\usepackage{epsfig}
\usepackage{amsthm}
\usepackage{amssymb}
\usepackage{amsmath}
\usepackage{amscd}
\usepackage{color}
\usepackage[T1]{fontenc}
\usepackage{upgreek}
\usepackage[font=scriptsize]{caption}
\usepackage{wrapfig}
\usepackage{graphicx}
\usepackage{subfigure}

\newcommand{\bg}{\begin{equation}}
\newcommand{\ed}{\end{equation}}
\newcommand{\bga}{\begin{eqnarray}}
\newcommand{\eda}{\end{eqnarray}}
\newcommand{\pf}{\textbf{Proof:\ }}

\def\cbdu{\par{\raggedleft$\Box$\par}}

\newtheorem {Theorem}  {Theorem}

\numberwithin{Theorem}{section}

\newtheorem {Lemma}[Theorem]  {Lemma}

\theoremstyle{definition}

\theoremstyle{remark}
\newtheorem{Remark}[Theorem]{\bf Remark}

\expandafter\chardef\csname pre amssym.def
at\endcsname=\the\catcode`\@ \catcode`\@=11
\def\undefine#1{\let#1\undefined}
\def\newsymbol#1#2#3#4#5{\let\next@\relax
 \ifnum#2=\@ne\let\next@\msafam@\else
 \ifnum#2=\tw@\let\next@\msbfam@\fi\fi
 \mathchardef#1="#3\next@#4#5}
\def\mathhexbox@#1#2#3{\relax
 \ifmmode\mathpalette{}{\m@th\mathchar"#1#2#3}%
 \else\leavevmode\hbox{$\m@th\mathchar"#1#2#3$}\fi}
\def\hexnumber@#1{\ifcase#1 0\or 1\or 2\or 3\or 4\or 5\or 6\or 7\or 8\or
 9\or A\or B\or C\or D\or E\or F\fi}

\font\teneufm=eufm10 \font\seveneufm=eufm7 \font\fiveeufm=eufm5
\newfam\eufmfam
\textfont\eufmfam=\teneufm \scriptfont\eufmfam=\seveneufm
\scriptscriptfont\eufmfam=\fiveeufm

\catcode`\@=\csname pre amssym.def at\endcsname

\newcounter{remark}
\def  \12  {{\frac{1}{2}}}

\def\build#1_#2^#3{\mathrel{\mathop{\kern 0pt#1}\limits_{#2}^{#3}}}

\numberwithin{equation}{section}

\begin{document}
%\currannalsline{0}{2006}

%\title[Local existence for Hall-MHD]{Local existence for the non-resistive Hall-magneto-hydrodynamics system in $\R^n$}
\title[Reduced models for EMHD]{Reduced models for electron magnetohydrodynamics: well-posedness and singularity formation}

%\author{hello}

\author [Mimi Dai]{Mimi Dai}

\address{Department of Mathematics, Statistics and Computer Science, University of Illinois at Chicago, Chicago, IL 60607, USA}
\email{mdai@uic.edu}

\thanks{M. Dai is partially supported by the NSF grants DMS--1815069 and DMS--2009422, and the von Neumann Fellowship at the Institute for Advanced Study. }

\begin{abstract}

We propose some one-dimensional reduced models for the three-dimensional electron magnetohydrodynamics which involves a highly nonlinear Hall term with intricate structure. The models contain nonlocal nonlinear terms which are more singular than that of the one-dimensional models for the Euler equation and the surface quasi-geostrophic equation. Local well-posedness is obtained in certain circumstances. In contrast, for a model with nonlocal transport term, we show that finite time singularity may occur for some initial data.

\bigskip

KEY WORDS: magnetohydrodynamics; reduced models; well-posedness, singularity formation.

\hspace{0.02cm}CLASSIFICATION CODE: 35Q35, 65M70, 76D03, 76W05, 76B03.
\end{abstract}

\maketitle

\section{Introduction}
In order to capture the rapid magnetic reconnection phenomena in plasma physics, the authors of \cite{ADFL} rigorously derived the following incompressible magnetohydrodynamics (MHD) model with Hall effect 
\begin{equation}\label{mhd}
\begin{split}
u_t+(u\cdot\nabla) u-(B\cdot\nabla) B+\nabla P=&\ \nu\Delta u, \\
B_t+(u\cdot\nabla) B-(B\cdot\nabla) u 
+\nabla\times ((\nabla\times B)\times B)=&\ \mu\Delta B, \\
\nabla\cdot u=0, \ \ \nabla\cdot B=&\ 0,
\end{split}
\end{equation}
considered on the space time domain $\Omega\times [0,\infty)$ with $\Omega\subset \mathbb R^3$. In system (\ref{mhd}), the unknowns are the velocity field $u$, the magnetic field $B$ and the scalar pressure function $P$. The constant coefficients $\nu$ and $\mu$ denote respectively the kinetic viscosity and magnetic resistivity. The Hall term $\nabla\times ((\nabla\times B)\times B)$ appears to be more singular than other nonlinear terms in the system. Without it, system (\ref{mhd}) reduces to the classical MHD which shares many similarities with the Navier-Stokes equation (NSE). With the presence of the Hall term, (\ref{mhd}) no longer has a natural scaling as the classical MHD does. Instead, the Hall term introduces new scaling and geometry properties as well as new challenges. It seems desirable to understand the Hall term as a single target. As such, one can consider the electron MHD system
\begin{equation}\label{emhd}
\begin{split}
B_t+ \nabla\times ((\nabla\times B)\times B)=&\ \mu\Delta B, \\
\nabla\cdot B=&\ 0
\end{split}
\end{equation}
which is a special case of (\ref{mhd}) with $u\equiv 0$, that is the magnetohydrodynamics with static background ion flow. 
Note that the first equation of (\ref{emhd}) is quasi-linear. In general there are limited analytic tools to treat quasi-linear partial differential equations. Remarkably, the authors of \cite{JO} established various ill-posedness results for (\ref{emhd}) with $\mu=0$ by exploring the quasi-linear dispersive character mechanism. On the other hand, replacing $\Delta B$ by a hyerviscosity term $-(-\Delta)^{\alpha}$ with $\alpha>1$ in \cite{Dai22}, the author showed global existence of weak solutions almost surely for rough initial data. It was achieved through an appropriate randomization of the initial data.  Nevertheless, we encounter essential difficulties to lower down the value of $\alpha$ to be 1 regarding the almost sure global existence of weak solutions.  To further understand the intricate structure of the Hall term, we will propose and study some nonlocal nonlinear one-dimensional (1D) models for (\ref{emhd}). This is inspired by the vast investigation of 1D nonlocal models for the Euler equation, NSE and surface quasi-geostrophic (SQG) equation in the literature. 

\medskip

\subsection{1D models for fluid equations} 
As an attempt to understand the vorticity form of the Euler equation, a 1D model with nonlocal nonlinear structure was introduced by Constantin, Lax and Majda \cite{CLM}. This model was later generalized by De Gregorio \cite{DeG1, DeG2} and studied by many authors, for instance, see \cite{Ch1, Ch2, Ch3, CHH, CHK, CKY, EGM, EJ, LH, OSW, Sak}. Regarding these 1D toy models for the Euler equation and NSE, one important issue is the competition and balance between the stretching effect and the transport effect. With the efforts of the aforementioned authors, this question has been understood very well. Various well-posendess results were obtained in certain contexts, while solutions with finite time singularities were constructed in contrast scenarios. 

It is known that the 2D SQG has strikingly similar features with the 3D Euler equation. In this vein of gaining insights from simplified models, 1D reduced models for the SQG were also proposed and investigated \cite{BLM, CasC, CCCF, CCF1, CCF2, Do, Dong, LR1, LR2, Mor, SV}. Both the problems of well-posedness and finite time blowup have been studied in depth. Readers with interest are referred to the works above.

\medskip

\subsection{1D models for electron MHD} 
\label{sec-1dmhd}
% The toy models for the Hall MHD will be studied in a separate paper.
Following the footprint of the study on 1D models for the fluid equations, we propose a family of nonlocal nonlinear simplified models for the electron MHD (\ref{emhd}). 

Denote the current density by $J=\nabla\times B$. Note that $\nabla\cdot B=0$ and $\nabla\cdot J=0$. 
It follows from the Biot-Savart law that 
\begin{equation}\notag%\label{BS}
B=\nabla\times (-\Delta)^{-1} J.
\end{equation}
Equation (\ref{emhd}) with generalized dissipation can be written as 
\begin{equation}\label{emhd-1}
B_t+(B\cdot\nabla) J -(J\cdot\nabla) B+\mu \Lambda^\alpha B=0, \ \ \ \Lambda=(-\Delta)^{\frac12}
\end{equation}
with $\alpha\geq 0$.
We suggest a two-parameter family of 1D toy models to approximate (\ref{emhd-1}), 
%on $[-\pi, \pi]$,
\begin{equation}\label{emhd-1d}
\begin{split}
B_t+aB J_x+bJ B_x+\mu\Lambda^\alpha B=&\ 0,  \ \ \ \ a, b\in \mathbb R\\
B_x=&\ \mathcal H J
\end{split}
\end{equation}
with $\Lambda=\mathcal H\partial_x$ in one spatial dimension. The second equation of (\ref{emhd-1d}) is an 1D analog of the Biot-Savart law. 
We study (\ref{emhd-1d}) either on the periodic interval $\mathbb S^1=[-\pi, \pi]$ or on $\mathbb R$. The Hilbert transform $\mathcal H$ on $\mathbb S^1$ can be defined as
\begin{equation}\label{H2}
\mathcal H f(x)=\frac{1}{2\pi} P.V. \int_{-\pi}^{\pi}f(y)\cot\left( \frac{x-y}{2}\right)\, dy
\end{equation}
and on $\mathbb R$ as
\begin{equation}\label{H1}
\mathcal H f=\frac{1}{\pi} P.V. \int_{-\infty}^{\infty}\frac{f(y)}{x-y}\, dy.
\end{equation}
To uniquely determine $B$ from $J$, we make the choice of Gauge by taking zero-mean value
\begin{equation}\label{gauge1}
\int B(t,x)\,dx =0.
\end{equation}
%or fixing one point
%\begin{equation}\label{gauge2}
%p(t,0) = m(t, 0)=0.
%\end{equation}
%{\color{red} }
%Notice that 
%\[u=\frac12(p+m), \ \ \ B=\frac12(p-m),\]
%and hence
%\begin{equation}\label{ub-pm}
%\begin{split}
%u_x=\frac12(p_x+m_x)=& \frac12 H(\Omega+\omega), \\
%B_x=\frac12(p_x-m_x)=& \frac12 H(\Omega-\omega).
%\end{split}
%\end{equation}
Note that the mean value of $B$ is invariant for system (\ref{emhd-1d}) with any $a, b\in \mathbb R$. Indeed, we have  formally
\begin{equation}\notag
\begin{split}
\frac{d}{dt} \int B(t,x)\, dx=&\int \left(-aBJ_x-bJ B_x-\mu\Lambda^\alpha B\right)\,dx\\
=&\int (a-b)B_x J\,dx\\
=&\ (a-b)\int J \mathcal H J\,dx\\
=&\ 0
\end{split}
\end{equation}
where we have used the integration by parts and the skew symmetry property of the Hilbert transform. 

System (\ref{emhd-1d}) has the natural scaling, if $B(x,t)$ is a solution with initial data $B_0(x)$ then $B_\lambda(x,t):=\lambda^{\alpha-2} B(\lambda x, \lambda^\alpha t)$ with an arbitrary rescaling parameter $\lambda$ is a solution as well with initial data $\lambda^{\alpha-2} B_0(\lambda x)$. Under such scaling the Sobolev space $\dot H^{\frac52-\alpha}$ is scaling invariant, also called critical. Note that $L^2$ is critical when $\alpha=\frac52$. By convention, (\ref{emhd-1d}) is referred to be supercritical, critical and subcritical respectively for $\alpha<\frac52$, $\alpha=\frac52$ and $\alpha>\frac52$. 

In order to analyze the quadratic nonlinear terms in (\ref{emhd-1d}) separately, we also consider the model with $a=1$ and $b=0$
\begin{equation}\label{e1d2}
B_t+BJ_x+\mu\Lambda^\alpha B= 0, \ \ \ 
B_x= \mathcal H J
\end{equation}
and the one with $a=0$ and $b=1$
\begin{equation}\label{e1d3}
B_t+JB_x+\mu\Lambda^\alpha B= 0, \ \ \ 
B_x= \mathcal H J.
\end{equation}
%{\color{blue} We should study well-posedness of (\ref{e1d2}) and (\ref{e1d3}) separately. The well-posedness is valid on both $\mathbb R$ and $\mathbb S$. The proof is only on $\mathbb S$.}
%We can also consider the model in conservative form,
%\begin{equation}\label{e1d4}
%\begin{split}
%B_t+(BJ)_x+\mu\Lambda^\alpha B=&\ 0,\\
%B_x=&\ \mathcal H J.
%\end{split}
%\end{equation}
%Question: what is the benefit of being in conservative form? 
%Note that (\ref{emhd-1d}), (\ref{e1d2}), (\ref{e1d3}) and (\ref{e1d4}) exhibit even symmetry in the sense that $B(t, -x)$ is a solution provided $B(t,x)$ is a solution. Due to the conservation of $\int B \, dx$, it is appropriate to study solutions of (\ref{emhd-1d}), (\ref{e1d2}) and (\ref{e1d3}) in the class of zero-mean value functions. 
Both systems (\ref{e1d2}) and (\ref{e1d3}) inherit the same scaling property of (\ref{emhd-1d}). One obvious obstacle in analyzing (\ref{e1d2}) and (\ref{e1d3}) is that the flux terms 
\[\int BJ_xB\, dx \ \ \ \mbox{and} \ \ \ \int JB_x B\,dx\]
do not vanish for smooth function $B$ in general. Such lack of symmetry and cancellation makes it hard to obtain the a priori estimate for $B$ in $L^2$. It also causes difficulties in establishing estimates in spaces with higher order regularity. On the other hand, the presence of the dissipative term $\Lambda^\alpha$ which has a typical regularizing effect is beneficial to establish well-posedness. We notice that the nonlinear term $BJ_x$ is more singular than $JB_x$. Intuitively, one expects that solutions of (\ref{e1d3}) behave more regularly than that of (\ref{e1d2}). Indeed, we will show local well-posedness for (\ref{e1d2}), and in general for (\ref{emhd-1d}), with $\alpha>2$ (including a portion of the supercritical regime); while we are able to obtain local well-posedness for (\ref{e1d3}) with $\alpha>1$. The situation for $0<\alpha\leq 1$ is more challenging, in which case (\ref{e1d3}) is quasi-linear. Moreover, compared to the nonlinear structure of Burgers' equation or the 1D nonlocal transport equation studied in \cite{Dong}, $JB_x$ is more singular since the "velocity" term $J$ involves derivative of $B$. Combined with the fact of lack of cancellation symmetry, we find it difficult to establish local well-posedness for (\ref{e1d3}) with $0<\alpha\leq 1$. Note that the loss of derivative in (\ref{e1d3}) is one degree. The existence of local analytic solutions is expected. 

When $\mu=0$, (\ref{e1d3}) is a transport equation. In particular, the maximum principle holds in this situation. When $\mu>0$ and $0< \alpha<2$, we can show that the maximum principle remains valid for (\ref{e1d3}). With the help of maximum principle, we are able to show that for some initial data, solutions of (\ref{e1d3}) with $\mu=0$ can develop singularities in finite time. 

In general belief, global well-posedness may be obtained for critical and subcritical systems, with $\alpha\geq \frac52$ in our context. However, the aforementioned lack of cancellation symmetry prevents us to achieve the global well-posedness by either using standard energy method or the approach of continuity of moduli. 
In another direction, since $BJ_x$ is more singular than $JB_x$, we conjecture that there exists initial data such that solutions of (\ref{e1d2}) with appropriate $\alpha$ blows up in finite time. These questions will be addressed in future work.

\medskip

\subsection{Main results}
We first state the local well-posedness results for (\ref{emhd-1d}) (valid for (\ref{e1d2}) as well) and (\ref{e1d3}). Although stated for $\Omega=\mathbb S^1$, they are also true for $\Omega=\mathbb R$ with slightly modified proof. 

For general $a,b\in\mathbb R$, we show local well-posedness of (\ref{emhd-1d}) with $\alpha\in(2, \frac52]$ in the critical space $H^{\frac52-\alpha}(\mathbb S^1)$ and some smoothing estimates. 
\begin{Theorem}\label{thm-local}
Let $a, b\in \mathbb R$, $\alpha\in(2, \frac52]$ and $B_0\in  H^{\frac52-\alpha}(\mathbb S^1)$. There exists a time $T>0$ depending on $\|B_{0}\|_{H^{\frac52-\alpha}}$ such that there exists a unique solution $B(t,x)$ to (\ref{emhd-1d}) with initial data $B(0,x)=B_0$ on $[0,T)$ satisfying
\[B\in C\left([0,T);   H^{\frac52-\alpha}(\mathbb S^1)\right)\cap L^2\left([0,T);  H^{\frac52-\frac{\alpha}2}(\mathbb S^1)\right)\]
and 
\begin{equation}\label{est-b1}
\sup_{0<t<T} t^{\frac{\beta}{\alpha}} \|B(t)\|_{ H^{\frac52-\alpha+\beta}(\mathbb S^1)}<\infty \ \ \forall \ \ \beta\geq 0,
\end{equation}
\begin{equation}\label{est-b2}
\lim_{t\to 0} t^{\frac{\beta}{\alpha}} \|B(t)\|_{H^{\frac52-\alpha+\beta}(\mathbb S^1)}=0 \ \ \forall \ \ \beta> 0.
\end{equation}
\end{Theorem}

%{\color{blue}This may be improved to $\alpha\in(1, \frac52]$ by using a different integration by parts trick for the term $BJ_x$. }

\medskip
In the subcritical regime, we show
\begin{Theorem}\label{thm-local2}
Let $a, b\in \mathbb R$, $\alpha\in(\frac52, \infty)$ and $B_0\in  L^2(\mathbb S^1)$. There exists a time $T$ which depends on $\|B_{0}\|_{L^2}$ such that there exists a unique solution $B(t,x)$ to (\ref{emhd-1d}) with initial data $B(0,x)=B_0$ on $[0,T)$ satisfying
\[B\in C\left([0,T);   L^2(\mathbb S^1)\right)\cap L^2\left([0,T);  H^{\frac{\alpha}2}(\mathbb S^1)\right)\]
and 
\begin{equation}\label{est-b3}
\sup_{0<t<T} t^{\frac{\beta}{\alpha}} \|B(t)\|_{ H^{\beta}(\mathbb S^1)}<\infty \ \ \forall \ \ \beta\geq 0,
\end{equation}
\begin{equation}\label{est-b4}
\lim_{t\to 0} t^{\frac{\beta}{\alpha}} \|B(t)\|_{H^{\beta}(\mathbb S^1)}=0 \ \ \forall \ \ \beta> 0.
\end{equation}
\end{Theorem}

\medskip

Without the more singular term $BJ_x$, we are able to show local well-posedness for (\ref{e1d3}) with smaller value of $\alpha$. Namely,
\begin{Theorem}\label{thm-local1}
Let $\alpha\in(1, \frac52]$ and $B_0\in  H^{\frac52-\alpha}(\mathbb S^1)$. The statements of Theorem \ref{thm-local} are true for model (\ref{e1d3}).
\end{Theorem}

For small $\alpha$, say $\alpha\leq 1$, we prove the local existence of analytic solutions to (\ref{e1d3}). 

\begin{Theorem}\label{thm-analytic}
Let $\alpha\in[0, 1]$ and $\mu\geq 0$. Assume the initial data $B_0\in L^2(\mathbb R)\cap C^2(\mathbb R)$ vanishes at infinity. There exists a time $T>0$ such that equation (\ref{e1d3}) has a real analytic solution on $(0,T)$ with $B(x,0)=B_0(x)$.
\end{Theorem}

In the absence of $BJ_x$, with the help of maximum principle, we show finite time singularity occurs for (\ref{e1d3}) with $\mu=0$ and a class of initial data. 
%There might be some hope to show local well-posedness with the help of maximum principle. 

\begin{Theorem}\label{thm-s1}
Assume the initial profile $B_0$ is odd, smooth, nonnegative on $[0, \infty)$, compactly supported on $[-1,1]$ with $\max_{x\in\mathbb R}B_0(x)=1$ and satisfies $\partial_{xx}B_0\leq 0$ on $[0, \infty)$.  Let $B(x,t)$ be the solution of (\ref{e1d3}) with $\mu=0$ and the initial data $B_0$ obtained in Theorem \ref{thm-analytic}. Then the norm $\|B_x(t)\|_{L^\infty}$ blows up at a finite time.
%Let $B(x,t)$ be a solution to (\ref{eq-s1}) with $\mu=0$ and the initial data $B_0$. Then the norm $\|B_x\|_{L^\infty}$ blows up in finite time.
\end{Theorem}

%\begin{Theorem}\label{thm-s2}
%Let $0<\delta<1$ and $0\leq \alpha<\frac32-\frac{\delta}{2}$. Let $B_0$ be the profile satisfying the conditions in Theorem \ref{thm-s1}. In addition, we assume 
%\[\int_0^\infty \frac{B_0(x)}{x^{\delta}}\, dx>C\]
%for a constant $C>0$. The conclusion of Theorem \ref{thm-s1} holds for (\ref{e1d3}) with $\mu>0$ and the initial data $B_0$.
%\end{Theorem}

%{\color{blue} For real analytic initial data, one might be able to show local existence for the model by using the Cauchy-Kovalevskaya theorem. }

%{\color{blue} Fix the gap in the proof of the singularity formation.}

Note if $B$ satisfies (\ref{e1d3}), $-B$ satisfies 
\begin{equation}\label{e1d4}
B_t-B_xJ+\mu \Lambda^\alpha B=0,  \ \ B_x=\mathcal H J.
\end{equation} 
Hence the finite time singularity scenario in Theorems \ref{thm-s1} holds for (\ref{e1d4}) as well.

% It is not clear how to show  global well-posedness for (\ref{e1d3}) with $\alpha\geq \frac52$.
%Namely, we aim to show:
%\begin{Theorem}\label{thm-global1}
%Let $\alpha\geq \frac52$. Assume $B_0\in H^{\frac{\alpha}{2}}(\mathbb S^1)$. Then there exists a unique solution $B(t)$ of (\ref{e1d3}) with initial data $B_0$ in $L^\infty([0,\infty); H^{\frac{\alpha}{2}}(\mathbb S^1))\cap L^2([0,\infty); H^{\alpha}(\mathbb S^1))$.
%\end{Theorem}

%\begin{Theorem}\label{thm-global2}
%Let $\alpha\geq \frac52$. Assume $B_0\in H^{\frac{\alpha}{2}}(\mathbb R)$ and $\|B_0\|_{L^\infty}<\mu$. Then there exists a unique solution $B(t)$ of (\ref{e1d3}) with initial data $B_0$ in $L^\infty([0,\infty); H^{\frac{\alpha}{2}}(\mathbb R))\cap L^2([0,\infty); H^{\alpha}(\mathbb R))$.
%\end{Theorem}

%\begin{Theorem}\label{thm-global3}
%Let $\alpha\geq 4$. Assume $B_0\in H^{\frac{\alpha}{2}}(\mathbb S^1)$. Then there exists a unique solution $B(t)$ of (\ref{emhd-1d}) with initial data $B_0$ in $L^\infty([0,\infty); H^{\frac{\alpha}{2}}(\mathbb S^1))\cap L^2([0,\infty); H^{\alpha}(\mathbb S^1))$.
%\end{Theorem}

\bigskip

\section{Notations and preliminaries}

\subsection{Notations}
We denote $C$ by a general constant, $C(\delta)$ a constant depending on the parameter $\delta$ (it also applies to other parameters), all of which may vary from line to line. We adapt the simplified notation $\lesssim$ for $\leq$ up to a multiplication of constant.

\medskip

\subsection{Littlewood-Paley theory} We recall the notations for Littlewood-Paley theory on $\mathbb T^n=[-L, L]^n$ briefly. Denote $\lambda_q=2^q/L$ for any integer $q$. Let $\chi\in C_0^\infty(\mathbb R^n)$ be a nonnegative radial function defined as
\begin{equation}\notag
\chi(\xi)=
\begin{cases}
1, \ \ \ \mbox{for} \ \ |\xi|\leq \frac34\\
0, \ \ \ \mbox{for} \ \ |\xi|\geq 1.
\end{cases}
\end{equation}
Take $\varphi(\xi)=\chi(\frac{\xi}{2})-\chi(\xi)$ and 
\begin{equation}\notag
\varphi_q(\xi)=
\begin{cases}
\varphi(\lambda_q^{-1}\xi), \ \ \ \mbox{for} \ \ q\geq0\\
\chi(\xi), \ \ \ \ \ \ \ \  \mbox{for} \ \ q= -1. 
\end{cases}
\end{equation}
For a tempered distribution vector field $u$ on $\mathbb T^n$ and $q\geq -1$, we define the $q$-th Littlewood-Paley projection of $u$ 
\[u_q(x)=\Delta_q u(x):= \sum_{k\in \mathbb Z^n} \widehat u(k)\varphi_q(k) e^{i\frac{2\pi}{L} k\cdot x}\]
with $\widehat u(k)$ being the Fourier coefficient of $u$.  We have the decomposition 
\[u=\sum_{q=-1}^\infty u_q\]
in the distributional sense. To simplify notations, we use
\[u_{\leq Q}=\sum_{q\geq-1}^Q u_q, \ \ \ \ \ \ \widetilde {\Delta_q} u=u_{q-1}+u_q+u_{q+1}.\]
It is worth to mention that the norm $\|u\|_{H^s(\mathbb T^n)}$ is equivalent to 
\[\left(\sum_{q=-1}^\infty \lambda_q^{2s} \|u_q\|_{L^2}^2\right)^{\frac12}.\]

\medskip

\subsection{Auxiliary estimates}
\label{sec-basic}
We recollect some standard estimates regarding the operator $e^{-\mu t\Lambda^{\alpha}}$ in the following.
\begin{Lemma}\label{le-m1} \cite{Miu}
Let $\alpha>0$ and $f\in L^2$. There exist constants $C_1(\alpha), C_2(\alpha)>0$ such that 
\begin{equation}\notag
e^{-C_1(\alpha)\mu \lambda_q^{\alpha}t} \|f_q\|_{L^2}\leq \left\|e^{-\mu t\Lambda^{\alpha}}* f_q \right\|_{L^2}
\leq e^{-C_2(\alpha)\mu \lambda_q^{\alpha}t} \|f_q\|_{L^2}.
\end{equation}
\end{Lemma}

\begin{Lemma}\label{le-m2} \cite{Miu}
Let $\alpha>0$, $s\geq 0$ and $f\in L^2$. There exists a constant $C(\alpha, s)>0$ such that 
\begin{equation}\notag
\begin{split}
\sup_{t\in(0,\infty)}t^{\frac{s}{\alpha}} \left\|e^{-\mu t\Lambda^{\alpha}}* f \right\|_{H^s}\leq&\ C(\alpha, s) \|f\|_{L^2},\\
\lim_{t\to 0}t^{\frac{s}{\alpha}} \left\|e^{-\mu t\Lambda^{\alpha}}* f \right\|_{H^s}=&\ 0.
\end{split}
\end{equation}
In addition if $s\in[0,\frac{\alpha}{2}]$, we have
\begin{equation}\notag
 \left\|e^{-\mu t\Lambda^{\alpha}}* f \right\|_{L^{\frac{\alpha}{s}}_tH^s_x}\leq C(\alpha, s) \|f\|_{L^2}.
\end{equation}

\end{Lemma}

\medskip

Denote the commutator
\[[f, \Delta_q] g=f g_q-\Delta_q(fg). \]
We have the following commutator estimate.
\begin{Lemma}\label{le-d1} \cite{Dong}
Let $f$ and $g$ be functions on $\mathbb T^n$. 
%(i) Assume $m\geq 0$, $\frac{n}{2}\leq s_1<1+\frac{n}{2}$, $s_2<\frac{n}{2}$ and $m+s_1+s_2>0$. Assume $f\in  H^{s_1}\cap  H^{m+s_1}$ and $g\in  H^{s_2}\cap  H^{m+s_2}$. There exist a constant $C=C(m, n, s_1, s_2)>0$ and a sequence $\{c_q\}\in l^2$ with $\|c_q\|_{l^2}\leq 1$ such that for any $q\geq 0$
%\begin{equation}\notag
%\left\|[f, \Delta_q] g\right\|_{ H^m}\leq C c_q\lambda_q^{-(s_1+s_2-\frac{n}{2})} \left(\|f\|_{H^{m+s_1}}\|g\|_{H^{s_2}}+\|f\|_{H^{s_1}}\|g\|_{ H^{m+s_2}}\right).
%\end{equation} 
 Let $m\geq 0$, $s_1<1+\frac{n}{2}$, $s_2<\frac{n}{2}$ and $m+s_1+s_2>0$. Assume $f\in  H^{s_1}\cap  H^{m+s_1}$ and $g\in  H^{s_2}\cap  H^{m+s_2}$. There exist a constant $C=C(m, n, s_1, s_2)>0$ and a sequence $\{c_q\}\in l^2$ with $\|c_q\|_{l^2}\leq 1$ such that for any $q\geq 0$  
\begin{equation}\notag
\left\|\widetilde {\Delta_q}[f, \Delta_q] g\right\|_{L^2}\leq C c_q\lambda_q^{-(m+s_1+s_2-\frac{n}{2})} \left(\|f\|_{H^{m+s_1}}\|g\|_{H^{s_2}}+\|f\|_{H^{s_1}}\|g\|_{H^{m+s_2}}\right).
\end{equation} 
\end{Lemma}

\medskip

The next two lemmas are also import in the estimates later. 
\begin{Lemma}\label{le-d2} \cite{Dong}
Let $f$ and $g$ be functions on $\mathbb T^n$. 
Let $m\geq 0$, $s_1<\frac{n}{2}$ and $s_2\in \mathbb R$. Assume $f\in H^{s_1}\cap H^{m+s_1}$ and $g\in H^{s_2}\cap  H^{m+s_2}$. There exist a constant $C=C(m, n, s_1, s_2)>0$ and a sequence $\{c_q\}\in l^2$ with $\|c_q\|_{l^2}\leq 1$ such that for any $q\geq 0$
\begin{equation}\notag
\left\|\widetilde {\Delta_q} (fg_q)\right\|_{L^2}\leq C c_q\lambda_q^{-(m+s_1+s_2-\frac{n}{2})} \left(\|f\|_{H^{m+s_1}}\|g\|_{H^{s_2}}+\|f\|_{H^{s_1}}\|g\|_{H^{m+s_2}}\right).
\end{equation} 

\end{Lemma}

\medskip

\begin{Lemma}\label{le-d3} \cite{RS}
Let $s_1, s_2<\frac12$ and $s_1+s_2>0$. If $f\in H^{s_1}$ and $g\in H^{s_2}$ we have
\[\|fg\|_{H^{s_1+s_2-\frac{n}{2}}}\leq C \|f\|_{H^{s_1}} \|g\|_{H^{s_2}}\]
for a constant $C=C(n, s_1, s_2)>0$.
\end{Lemma}

\medskip

%\subsection{Properties of Hilbert transform}
%The Hilbert transform has the following simple properties
%\begin{equation}\notag
%\begin{split}
%\mathcal H(cf)=&\ c\mathcal Hf, \ \ \mbox{for a constant} \ \ c,\\
%\mathcal H\sin(kx)=&-\cos(kx),\ \ \ \mathcal H\cos(kx)=\sin(kx).
%\end{split}
%\end{equation}
%And more generally, we have
%\begin{equation}\notag
%\mathcal H\sin(kx+\theta)=-\cos(kx+\theta), \ \ \mathcal H\cos(kx+\theta)=\sin(kx+\theta).
%\end{equation}

%\begin{equation}\notag
%\begin{split}
%\mathcal H(\mathcal H f)=&-f,\\
%\mathcal H(fg)=&\ f\mathcal H g+g\mathcal H f+\mathcal H(\mathcal Hf\mathcal Hg),\\
%\mathcal H(f\mathcal H f)=&\ \frac12\left((\mathcal H f)^2-f^2\right), \\
%\mathcal H(e^{ikx})=&\ i \cdot\mathrm{sign}(k) e^{ikx}.
%\end{split}
%\end{equation}

%For any periodic function $f$, the mean value of its Hilbert transform is zero, that is
%\begin{equation}\label{Hm}
%\int \mathcal Hf\, dx=0.
%\end{equation}

We conclude this section by an estimate of the Hilbert transform. 

\begin{Lemma}\cite{Zy}
The Hilbert transform $\mathcal H$ is a bounded linear operator from space $L^p$ to $L^p$ with $1<p<\infty$ and
\begin{equation}\label{HLp}
\|\mathcal Hf\|_{L^p}\leq C(p) \|f\|_{L^p}
\end{equation}
for a constant $C(p)>0$.
\end{Lemma}

\bigskip

\section{Local well-posedness in the supercritical and critical regime}
\label{sec-local}

In this section, we prove the well-posedness results stated in Theorems \ref{thm-local} and \ref{thm-local1} for $\alpha\leq \frac52$. This is achieved in three steps: we first establish the a priori estimates; we then solve the approximating system, see (\ref{sys-app}), to obtain a sequence of approximating solutions satisfying the a priori estimates; in the end, we show the sequence is a Cauchy sequence in appropriate functional spaces by employing the a priori estimates and hence converges to a solution. The main difficulties arise in the first step. As explained in Subsection \ref{sec-1dmhd}, the fluxes
\[\int BJ_x B\, dx \ \ \ \mbox{and} \ \ \ \int JB_xB\, dx\] 
do not vanish due to lack of cancellation. It prevents us obtaining an estimate in $L^2$ immediately, opposed to the original PDE system. A similar situation happens for the 1D transport nonlocal equation studied in \cite{Dong}. The author overcome it by utilizing a backward bootstrap argument, in which higher order a priori estimates were first obtained and then lower order estimates were derived from the higher order ones. We adapt this type of backward bootstrap strategy for our models. However, the nonlinear terms $BJ_x$ and $JB_x$ in our models are more singular than that of the model from \cite{Dong}. Moreover, the term $BJ_x$ is not a transport term in (\ref{emhd-1d}). These features cause more difficulties in establishing the a priori estimates. They are the main reasons we are not able to reach the quasi-linear regime in terms of well-posedness, that is $\alpha\leq 2$ for (\ref{emhd-1d}) and $\alpha\leq 1$ for (\ref{e1d3}).

%\medskip

\subsection{A priori estimates in $H^s$}
\label{sec-priori}

Acting the projection $\Delta_q$ on equation (\ref{emhd-1d}) yields
\begin{equation}\notag
\partial_t B_q+a\Delta_q(BJ_x)+b\Delta_q(JB_x)+\mu\Lambda^\alpha B_q=0.
\end{equation}
Applying the commutation notation we have
\begin{equation}\notag
\partial_t B_q+aBJ_{x,q}+bJB_{x,q}+\mu\Lambda^\alpha B_q=a[B, \Delta_q]J_x-b[J, \Delta_q]B_x.
\end{equation}
Multiplying the equation above by $B_q$ and integrating over $\mathbb S^1$ gives
\begin{equation}\notag
\begin{split}
&\frac12\frac{d}{dt}\|B_q\|_{L^2}^2+\mu\lambda_q^\alpha \|B_q\|_{L^2}^2\\
\leq &\ a\int_{\mathbb S^1} [B, \Delta_q]J_x B_q\, dx+b\int_{\mathbb S^1} [J, \Delta_q]B_x B_q\, dx\\
&-a\int_{\mathbb S^1} BJ_{x,q} B_q\, dx-b\int_{\mathbb S^1} JB_{x,q} B_q\, dx\\
=&\ a \int_{\mathbb S^1} \widetilde{\Delta_q}\left([B, \Delta_q]J_x\right) B_q\, dx+b\int_{\mathbb S^1} \widetilde{\Delta_q}\left([J, \Delta_q]B_x\right) B_q\, dx\\
&-a\int_{\mathbb S^1} \widetilde{\Delta_q}\left(BJ_{x,q}\right) B_q\, dx-b\int_{\mathbb S^1} \widetilde{\Delta_q}\left(JB_{x,q}\right) B_q\, dx.
\end{split}
\end{equation}
Applying H\"older's inequality we obtain
\begin{equation}\label{energy1}
\begin{split}
&\frac{d}{dt}\|B_q\|_{L^2}+\mu\lambda_q^\alpha \|B_q\|_{L^2}\\
\leq &\ a \left\|\widetilde{\Delta_q}\left([B, \Delta_q]J_x\right)\right\|_{L^2} + b\left\| \widetilde{\Delta_q}\left([J, \Delta_q]B_x\right) \right\|_{L^2} \\
&+ a\left\| \widetilde{\Delta_q}\left(BJ_{x,q}\right) \right\|_{L^2} + b\left\| \widetilde{\Delta_q}\left(JB_{x,q}\right) \right\|_{L^2} 
\end{split}
\end{equation}
and 
\begin{equation}\label{energy2}
\begin{split}
&\frac{d}{dt}\|B_q\|_{L^2}+\mu\lambda_q^\alpha \|B_q\|_{L^2}\\
\geq &- a\left\|\widetilde{\Delta_q}\left([B, \Delta_q]J_x\right)\right\|_{L^2} - b\left\| \widetilde{\Delta_q}\left([J, \Delta_q]B_x\right) \right\|_{L^2} \\
&-a \left\| \widetilde{\Delta_q}\left(BJ_{x,q}\right) \right\|_{L^2} - b \left\| \widetilde{\Delta_q}\left(JB_{x,q}\right) \right\|_{L^2}.
\end{split}
\end{equation}
In the following we show an upper bound from (\ref{energy1}). A lower bound can be established analogously from (\ref{energy2}). First of all, we have from (\ref{energy1}) for any $s\geq 0$
\begin{equation}\label{energy3}
\begin{split}
\lambda_q^s\|B_q(t)\|_{L^2}
\leq&\ \lambda_q^se^{-\mu\lambda_q^\alpha t} \|B_q(0)\|_{L^2}\\
&+ a\int_0^t \lambda_q^se^{-\mu\lambda_q^\alpha (t-\tau)}\left\|\widetilde{\Delta_q}\left([B, \Delta_q]J_x\right)(\tau)\right\|_{L^2} \, d\tau\\
&+ b \int_0^t \lambda_q^se^{-\mu\lambda_q^\alpha (t-\tau)}\left\| \widetilde{\Delta_q}\left([J, \Delta_q]B_x\right) (\tau)\right\|_{L^2}\, d\tau \\
&+ a \int_0^t \lambda_q^se^{-\mu\lambda_q^\alpha (t-\tau)}\left\| \widetilde{\Delta_q}\left(BJ_{x,q}\right)(\tau) \right\|_{L^2}\, d\tau\\
& + b\int_0^t \lambda_q^se^{-\mu\lambda_q^\alpha (t-\tau)} \left\| \widetilde{\Delta_q}\left(JB_{x,q}\right)(\tau) \right\|_{L^2} \, d\tau.
\end{split}
\end{equation}
It follows from Lemma \ref{le-m1} and Lemma \ref{le-m2} that for $\frac52-\alpha\leq s\leq \frac{\alpha}{2}$ 
%and $2\leq \alpha\leq \frac52$
\begin{equation}\label{term1}
\begin{split}
&\left\| \lambda_q^se^{-\mu\lambda_q^\alpha t} \|B_q(0)\|_{L^2}\right\|_{L^{\frac{\alpha}{s+\alpha-\frac52}}(0,T; l^2)}\\
\leq &\ \left\| e^{-\mu\lambda_q^\alpha t} B_q(0)\right\|_{L^{\frac{\alpha}{s+\alpha-\frac52}}(0,T; H^s)}\\
\leq&\ C(T) \|B(0)\|_{H^{\frac{5}{2} -\alpha}}
\end{split}
\end{equation}
where the constant $C(T)\to 0$ as $T\to 0$. We apply Lemma \ref{le-d1} with $m=0$, $s_1=s<\frac32$ and $s_2=s-2<\frac12$ to the second term of the right hand side of (\ref{energy3}) 
\begin{equation}\notag
\begin{split}
& \int_0^t \lambda_q^se^{-\mu\lambda_q^\alpha (t-\tau)}\left\|\widetilde{\Delta_q}\left([B, \Delta_q]J_x\right)(\tau)\right\|_{L^2} \, d\tau\\
\leq & \int_0^t c_q\lambda_q^{\frac52-s}e^{-\mu\lambda_q^\alpha (t-\tau)} \|B\|_{H^s}\|J_x\|_{H^{s-2}} \, d\tau\\
\leq & \int_0^t c_q(t-\tau)^{-\frac{1}{\alpha}(\frac52-s)} \|B\|_{H^s}^2 \, d\tau\\
\end{split}
\end{equation}
where we have used the fact $x^ae^{-x}\leq C$ for $a\geq 0$. Thus we apply Hardy-Littlewood-Sobolev inequality for $\frac52-\alpha< s<\frac32$
\begin{equation}\label{term2}
\begin{split}
&\left\| \int_0^t \lambda_q^se^{-\mu\lambda_q^\alpha (t-\tau)}\left\|\widetilde{\Delta_q}\left([B, \Delta_q]J_x\right)(\tau)\right\|_{L^2} \, d\tau\right\|_{L^{\frac{\alpha}{s+\alpha-\frac52}}(0,T; l^2)}\\
\lesssim & \left\| \int_0^t (t-\tau)^{-\frac{1}{\alpha}(\frac52-s)} \|B\|_{H^s}^2 \, d\tau\right\|_{L^{\frac{\alpha}{s+\alpha-\frac52}}(0,T)}\\
\lesssim & \left\|\|B\|_{H^s}^2\right\|_{L^{\frac{\alpha}{2(s+\alpha-\frac52)}}(0,T)}\\
\lesssim & \|B\|^2_{L^{\frac{\alpha}{(s+\alpha-\frac52)}}(0,T; H^s)}.
\end{split}
\end{equation}
Analogously, applying Lemma \ref{le-d1} with $m=0$, $s_1=s_2=s-1<\frac12$ to the third term on the right hand side of (\ref{energy3}), applying Lemma \ref{le-d2} with $m=0$, $s_1=s<\frac12$ and $s_2=s-2$ to the fourth term, and applying Lemma \ref{le-d2} with $m=0$, $s_1=s_2=s-1<\frac12$ to the fifth term, we claim similar estimates for these terms. 

Summarizing the analysis above, if $a=1$ and $b=0$, we have for $\frac52-\alpha<s< \frac{1}{2}$ and $\alpha> 2$
%$2\leq \alpha\leq \frac52$  
\begin{equation}\label{energy4}
\|B\|_{L^{\frac{\alpha}{(s+\alpha-\frac52)}}(0,T; H^s)}\leq C(T) \|B(0)\|_{H^{\frac{5}{2} -\alpha}}+C \|B\|^2_{L^{\frac{\alpha}{(s+\alpha-\frac52)}}(0,T; H^s)}.
\end{equation}
In the case of $a=0$ and $b=1$, (\ref{energy4}) holds for $\frac52-\alpha<s< \frac{3}{2}$ and $\alpha> 1$.

Now multiplying (\ref{energy1}) by $\lambda_q^{\frac52-\frac{\alpha}{2}}$ and applying Lemmas \ref{le-d1} and \ref{le-d2} the same way as above, we obtain
\begin{equation}\notag
\|B(t)\|_{H^{\frac52-\frac{\alpha}{2}}}\lesssim C(T) \|B(0)\|_{H^{\frac{5}{2} -\alpha}}+\int_0^t (t-\tau)^{-\frac{1}{\alpha}(5-\frac{\alpha}{2}-2s)} \|B(\tau)\|_{H^s}^2 \, d\tau
\end{equation}
for $\frac52-\alpha<s<\frac52-\frac{\alpha}{4}$. Taking $L^2$ norm in time and applying the Hardy-Littlewood inequality to the inequality above yields
\begin{equation}\label{energy5}
\|B\|_{L^{2}(0,T; H^{\frac52-\frac{\alpha}{2}})}\leq C(T) \|B(0)\|_{H^{\frac{5}{2} -\alpha}}+C \|B\|^2_{L^{\frac{\alpha}{(s+\alpha-\frac52)}}(0,T; H^s)}.
\end{equation}
Similarly we can show
\begin{equation}\label{energy6}
\|B\|_{L^{\infty}(0,T; H^{\frac52-\alpha})}\leq  \|B(0)\|_{H^{\frac{5}{2} -\alpha}}+C \|B\|^2_{L^{2}(0,T; H^{\frac52-\frac{\alpha}{2}})}.
\end{equation}
%{\color{blue}Add details to obtain the estimates...}

Note that we can establish a lower bound on $\|B\|_{L^{\infty}(0,T; H^{\frac52-\alpha})}$ by analogous analysis as above
\begin{equation}\label{energy-low}
\begin{split}
&\|B\|_{L^{\infty}(0,T; H^{\frac52-\alpha})}\\
\geq& \left(\sum_{q\geq0} \lambda_q^{2(\frac52-\alpha)} e^{-2\mu\lambda_q^\alpha t}\|B_q(0)\|^2_{L^2}\right)^{\frac12} -C \|B\|^2_{L^{2}(0,T; H^{\frac52-\frac{\alpha}{2}})}.
\end{split}
\end{equation}

\medskip

\subsection{Smoothing estimates (\ref{est-b1}) and (\ref{est-b2})}
\label{sec-smooth}
%Now we show (\ref{est-b1}) and (\ref{est-b2}). 
We first show the estimates (\ref{est-b1}) and (\ref{est-b2}) hold for $0\leq\beta\leq \frac{\alpha}{2}$ and then iterate the process to prove them for all $\beta>0$.  Taking $s=\frac{5}{2}-\alpha+\beta$ for $0\leq\beta\leq \frac{\alpha}{2}$ in (\ref{energy3}) and taking sum in $q$ yields
\begin{equation}\label{energy7}
\begin{split}
\|B(t)\|_{H^{\frac{5}{2}-\alpha+\beta}}
\leq&\ \sum_{q\geq0}\lambda_q^{\frac{5}{2}-\alpha+\beta}e^{-\mu\lambda_q^\alpha t} \|B_q(0)\|_{L^2}\\
&+ \sum_{q\geq0}\int_0^t \lambda_q^{\frac{5}{2}-\alpha+\beta}e^{-\mu\lambda_q^\alpha (t-\tau)}\left\|\widetilde{\Delta_q}\left([B, \Delta_q]J_x\right)(\tau)\right\|_{L^2} \, d\tau\\
&+  \sum_{q\geq0}\int_0^t \lambda_q^{\frac{5}{2}-\alpha+\beta}e^{-\mu\lambda_q^\alpha (t-\tau)}\left\| \widetilde{\Delta_q}\left([J, \Delta_q]B_x\right) (\tau)\right\|_{L^2}\, d\tau \\
&+ \sum_{q\geq0} \int_0^t \lambda_q^{\frac{5}{2}-\alpha+\beta}e^{-\mu\lambda_q^\alpha (t-\tau)}\left\| \widetilde{\Delta_q}\left(BJ_{x,q}\right)(\tau) \right\|_{L^2}\, d\tau\\
& + \sum_{q\geq0}\int_0^t \lambda_q^{\frac{5}{2}-\alpha+\beta}e^{-\mu\lambda_q^\alpha (t-\tau)} \left\| \widetilde{\Delta_q}\left(JB_{x,q}\right)(\tau) \right\|_{L^2} \, d\tau\\
=:&\ I_1+I_2+I_3+I_4+I_5.
\end{split}
\end{equation}
In view of the fact $x^ae^{-x}\leq C$ for $a\geq 0$, Lemma \ref{le-m1} and Lemma \ref{le-m2} again, we have
\begin{equation}\notag
I_1\leq C(T)\sum_{q\geq0}t^{-\frac{\beta}{\alpha}}\lambda_q^{\frac{5}{2}-\alpha}  \|B_q(0)\|_{L^2}
\leq C(T)t^{-\frac{\beta}{\alpha}} \|B(0)\|_{H^{\frac{5}{2}-\alpha} }.
\end{equation}
In order to estimate $I_2$, we apply Lemma \ref{le-d1} with $m=0$, $s_1=\frac52-\alpha+\beta<\frac32$ and $s_2=\frac12-\alpha+\beta<\frac12$ to infer
\begin{equation}\notag
\begin{split}
I_2\leq& \sum_{q\geq0}\int_0^t c_q\lambda_q^{\alpha-\beta}e^{-\mu\lambda_q^\alpha (t-\tau)}\|B(\tau)\|^2_{H^{\frac52-\alpha+\beta}} \, d\tau\\
\lesssim& \int_0^t (t-\tau)^{-1+\frac{\beta}{\alpha}}\|B(\tau)\|^2_{H^{\frac52-\alpha+\beta}} \, d\tau.
\end{split}
\end{equation}
Similarly we apply Lemma \ref{le-d1} with $m=0$, $s_1=s_2=\frac32-\alpha+\beta<\frac12$ to $I_3$, Lemma \ref{le-d2} with $m=0$, $s_1=\frac52-\alpha+\beta<\frac12$ and $s_2=\frac12-\alpha+\beta$ to $I_4$ and Lemma \ref{le-d2} with $m=0$, $s_1=s_2=\frac32-\alpha+\beta<\frac12$ to $I_5$ to obtain 
\begin{equation}\notag
I_3+I_4+I_5
\lesssim \int_0^t (t-\tau)^{-1+\frac{\beta}{\alpha}}\|B(\tau)\|^2_{H^{\frac52-\alpha+\beta}} \, d\tau.
\end{equation}
Combining the estimates above with (\ref{energy7}) we conclude that for $0\leq \beta<\alpha-2$ 
\begin{equation}\notag
\|B(t)\|_{H^{\frac{5}{2}-\alpha+\beta}}
\leq C(T)t^{-\frac{\beta}{\alpha}} \|B(0)\|_{H^{\frac{5}{2}-\alpha} }+C\int_0^t (t-\tau)^{-1+\frac{\beta}{\alpha}}\|B(\tau)\|^2_{H^{\frac52-\alpha+\beta}} \, d\tau
\end{equation}
immediately followed by 
\begin{equation}\notag
\begin{split}
&t^{\frac{\beta}{\alpha}}\|B(t)\|_{H^{\frac{5}{2}-\alpha+\beta}}\\
\leq&\ C(T) \|B(0)\|_{H^{\frac{5}{2}-\alpha} }+C\int_0^t t^{\frac{\beta}{\alpha}}(t-\tau)^{-1+\frac{\beta}{\alpha}}\|B(\tau)\|^2_{H^{\frac52-\alpha+\beta}} \, d\tau\\
\leq&\ C(T) \|B(0)\|_{H^{\frac{5}{2}-\alpha} }+C\left(\sup_{t\in(0,T)}t^{\frac{\beta}{\alpha}}\|B(t)\|_{H^{\frac{5}{2}-\alpha+\beta}}\right)^2\int_0^t t^{\frac{\beta}{\alpha}}(t-\tau)^{-1+\frac{\beta}{\alpha}}\tau^{-\frac{2\beta}{\alpha}} \, d\tau.
\end{split}
\end{equation}
Changing variable $\tau=t\tau'$ in the time integral we have
\begin{equation}\notag
\begin{split}
&\int_0^t t^{\frac{\beta}{\alpha}}(t-\tau)^{-1+\frac{\beta}{\alpha}}\tau^{-\frac{2\beta}{\alpha}} \, d\tau\\
=&\int_0^1 (1-\tau')^{-1+\frac{\beta}{\alpha}} (\tau')^{-1+(1-\frac{2\beta}{\alpha})}\, d\tau'\\
\lesssim&\ 1
\end{split}
\end{equation}
provided $\beta<\frac{\alpha}{2}$ and hence $1-\frac{2\beta}{\alpha}>0$. Therefore, for $0\leq\beta<\min\{\alpha-2, \frac{\alpha}{2}\}$, it follows
\begin{equation}\notag
\begin{split}
&\sup_{t\in(0,T)}t^{\frac{\beta}{\alpha}}\|B(t)\|_{H^{\frac{5}{2}-\alpha+\beta}}\\
\leq&\ C(T) \|B(0)\|_{H^{\frac{5}{2}-\alpha} }+CC(T)\left(\sup_{t\in(0,T)}t^{\frac{\beta}{\alpha}}\|B(t)\|_{H^{\frac{5}{2}-\alpha+\beta}}\right)^2
\end{split}
\end{equation}
where $C(T)\to 0$ as $T\to 0$. Thus for small enough $T>0$, estimates (\ref{est-b1}) and (\ref{est-b2}) hold. 

Now we assume $\beta\geq \min\{\alpha-2, \frac{\alpha}{2}\}$. We split $I_2$ as 
\[I_2= \sum_{q\geq 0}\int_0^{\frac{t}{2}} \cdot\cdot\cdot + \sum_{q\geq 0}\int_{\frac{t}{2}}^t \cdot\cdot\cdot =: I_{21} +I_{22}\]
and estimate $I_{21}$ and $I_{22}$ separately as follows. Applying Lemma \ref{le-d1} with $m=0$, $s_1=\frac52-\alpha+\gamma<\frac32$ and $s_2=\frac12-\alpha+\gamma<\frac12$ for $0<\gamma<\min\{\alpha-2, \frac{\alpha}{2}\}$ to $I_{21}$, we deduce
\begin{equation}\notag
\begin{split}
I_{21}\leq& \sum_{q\geq0}\int_0^{\frac{t}2}c_q\lambda_q^{\alpha+\beta-2\gamma}e^{-\mu \lambda_q^\alpha (t-\tau)} \|B(\tau)\|^2_{H^{\frac52-\alpha+\gamma}}\, d\tau\\
\lesssim& \int_0^{\frac{t}2} (t-\tau)^{-\frac{\alpha+\beta-2\gamma}{\alpha}} \|B(\tau)\|^2_{H^{\frac52-\alpha+\gamma}}\, d\tau.
\end{split}
\end{equation}
While for $I_{22}$ we apply Lemma \ref{le-d1} with $m>0$, $s_1=\frac52-\alpha+\gamma<\frac32$ and $s_2=\frac12-\alpha+\gamma<\frac12$ and obtain
\begin{equation}\notag
\begin{split}
I_{22}\leq& \sum_{q\geq0}\int_{\frac{t}2}^tc_q\lambda_q^{\alpha+\beta-2\gamma-m}e^{-\mu \lambda_q^\alpha (t-\tau)} \|B(\tau)\|_{H^{\frac52-\alpha+\gamma}}  \|B(\tau)\|_{H^{m+\frac52-\alpha+\gamma}}\, d\tau\\
\lesssim& \int_{\frac{t}2}^t (t-\tau)^{-\frac{\alpha+\beta-2\gamma-m}{\alpha}} \|B(\tau)\|_{H^{\frac52-\alpha+\gamma}}\|B(\tau)\|_{H^{m+\frac52-\alpha+\gamma}}\, d\tau.
\end{split}
\end{equation}
Similarly we split $I_j=I_{j1}+I_{j2}$ for $3\leq j\leq 5$. Analogous estimates as before shows that by applying \\
(i) Lemma \ref{le-d1}  to $I_{31}$ with $m=0$, $s_1=s_2=\frac32-\alpha+\gamma<\frac12$,\\
(ii) Lemma \ref{le-d2} to $I_{41}$ with $m=0$, $s_1=\frac52-\alpha+\gamma<\frac12$ and $s_2=\frac12-\alpha+\gamma$,\\
(iii) Lemma \ref{le-d2} to $I_{51}$ with $m=0$, $s_1=s_2=\frac32-\alpha+\gamma<\frac12$,\\
(iv) Lemma \ref{le-d1} to $I_{32}$ with $m>0$, $s_1=s_2=\frac32-\alpha+\gamma<\frac12$,\\
(v) Lemma \ref{le-d2} to $I_{42}$ with $m>0$, $s_1=\frac52-\alpha+\gamma<\frac12$ and $s_2=\frac12-\alpha+\gamma$,\\
(vi) Lemma \ref{le-d2} to $I_{52}$ with $m>0$, $s_1=s_2=\frac32-\alpha+\gamma<\frac12$,\\
we have 
\begin{equation}\notag
\begin{split}
I_{31}+I_{41}+I_{51}
\lesssim& \int_0^{\frac{t}2} (t-\tau)^{-\frac{\alpha+\beta-2\gamma}{\alpha}} \|B(\tau)\|^2_{H^{\frac52-\alpha+\gamma}}\, d\tau,\\
I_{32}+I_{42}+I_{52}
\lesssim& \int_{\frac{t}2}^t (t-\tau)^{-\frac{\alpha+\beta-2\gamma-m}{\alpha}} \|B(\tau)\|_{H^{\frac52-\alpha+\gamma}}\|B(\tau)\|_{H^{m+\frac52-\alpha+\gamma}}\, d\tau.
\end{split}
\end{equation}
Summarizing the analysis above we conclude that for $\beta\geq \min\{\frac52-\alpha, \frac{\alpha}{2}\}$, $0<\gamma< \min\{\frac52-\alpha, \frac{\alpha}{2}\}$ and $m>0$, 
\begin{equation}\notag
\begin{split}
&\|B(t)\|_{H^{\frac{5}{2}-\alpha+\beta}}\\
\leq&\ C(T)t^{-\frac{\beta}{\alpha}} \|B(0)\|_{H^{\frac{5}{2}-\alpha} }+C\int_0^{\frac{t}{2}} (t-\tau)^{-1+\frac{2\gamma-\beta}{\alpha}}\|B(\tau)\|^2_{H^{\frac52-\alpha+\gamma}} \, d\tau\\
&+C \int_{\frac{t}2}^t (t-\tau)^{-1+\frac{2\gamma+m-\beta}{\alpha}} \|B(\tau)\|_{H^{\frac52-\alpha+\gamma}}\|B(\tau)\|_{H^{m+\frac52-\alpha+\gamma}}\, d\tau.
\end{split}
\end{equation}
It implies 
\begin{equation}\label{est-high1}
\begin{split}
&\sup_{t\in(0,T)}t^{\frac{\beta}{\alpha}}\|B(t)\|_{H^{\frac{5}{2}-\alpha+\beta}}\\
%\leq&\ c(T) \|B(0)\|_{H^{\frac{5}{2}-\alpha} }+C\int_0^{\frac{t}{2}} t^{\frac{\beta}{\alpha}}(t-\tau)^{-1+\frac{2\gamma-\beta}{\alpha}}\|B(\tau)\|^2_{H^{\frac52-\alpha+\gamma}} \, d\tau\\
\leq&\ C(T) \|B(0)\|_{H^{\frac{5}{2}-\alpha} }\\
&+C\left(\sup_{t\in(0,T)}t^{\frac{\gamma}{\alpha}}\|B(t)\|_{H^{\frac{5}{2}-\alpha+\gamma}}\right)^2\int_0^{\frac{t}{2}} t^{\frac{\beta}{\alpha}}(t-\tau)^{-1+\frac{2\gamma-\beta}{\alpha}}\tau^{-\frac{2\gamma}{\alpha}} \, d\tau\\
&+C\left(\sup_{t\in(0,T)}t^{\frac{\gamma}{\alpha}}\|B(t)\|_{H^{\frac{5}{2}-\alpha+\gamma}}\right)
\left(\sup_{t\in(0,T)}t^{\frac{\gamma+m}{\alpha}}\|B(t)\|_{H^{m+\frac{5}{2}-\alpha+\gamma}}\right)\\
&\cdot\int_{\frac{t}{2}}^t t^{\frac{\beta}{\alpha}}(t-\tau)^{-1+\frac{2\gamma+m-\beta}{\alpha}}\tau^{-\frac{2\gamma+m}{\alpha}} \, d\tau.
\end{split}
\end{equation}
The first time integral is estimated by change of the variables
\begin{equation}\notag
\begin{split}
&\int_0^{\frac{t}{2}} t^{\frac{\beta}{\alpha}}(t-\tau)^{-1+\frac{2\gamma-\beta}{\alpha}}\tau^{-\frac{2\gamma}{\alpha}} \, d\tau\\
=&\int_0^{\frac{1}{2}} (1-\tau')^{-1+\frac{2\gamma-\beta}{\alpha}}(\tau')^{-\frac{2\gamma}{\alpha}} \, d\tau'\\
\lesssim& \int_0^{\frac{1}{2}}(\tau')^{-\frac{2\gamma}{\alpha}} \, d\tau'
\lesssim 1
\end{split}
\end{equation}
since $0<\frac{2\gamma}{\alpha}<1$. The second time integral satisfies
\begin{equation}\notag
\begin{split}
&\int_{\frac{t}{2}}^t t^{\frac{\beta}{\alpha}}(t-\tau)^{-1+\frac{2\gamma+m-\beta}{\alpha}}\tau^{-\frac{2\gamma+m}{\alpha}} \, d\tau\\
=& \int_{\frac{1}{2}}^1 (1-\tau')^{-1+\frac{2\gamma+m-\beta}{\alpha}}(\tau')^{-\frac{2\gamma+m}{\alpha}} \, d\tau'\\
\lesssim & \int_{\frac{1}{2}}^1 (1-\tau')^{-1+\frac{2\gamma+m-\beta}{\alpha}}\, d\tau' \lesssim 1
\end{split}
\end{equation}
where we require $2\gamma+m>\beta$ to obtain the last step. On the other hand, in view of previous estimates, we note for $\gamma+m<\min\{\beta, \frac52-\alpha, \frac{\alpha}{2}\}$
\begin{equation}\notag
\begin{split}
\sup_{t\in(0,T)}t^{\frac{\gamma}{\alpha}}\|B(t)\|_{H^{\frac{5}{2}-\alpha+\gamma}}<&\ \infty,\\
\sup_{t\in(0,T)}t^{\frac{\gamma+m}{\alpha}}\|B(t)\|_{H^{m+\frac{5}{2}-\alpha+\gamma}}<&\ \infty,\\
\lim_{t\to 0}t^{\frac{\gamma}{\alpha}}\|B(t)\|_{H^{\frac{5}{2}-\alpha+\gamma}}=&\ 0,\\
\lim_{t\to 0}t^{\frac{\gamma+m}{\alpha}}\|B(t)\|_{H^{m+\frac{5}{2}-\alpha+\gamma}}=&\ 0.
\end{split}
\end{equation}
Therefore, from (\ref{est-high1}) we conclude the estimates (\ref{est-b1}) and (\ref{est-b2}) for appropriate parameters $\gamma$, $m$ and 
\[\gamma+m<\beta<2\gamma+m, \ \ \ \beta\geq \min\left\{\frac52-\alpha, \frac{\alpha}{2}\right\}.\]
We can iterate this process to obtain (\ref{est-b1}) and (\ref{est-b2}) for all $\beta>0$.

\medskip

\subsection{A priori estimate in $L^2$} 
\label{sec-l2}

%Question: {\color{red} We are on torus, do we need this part below to show estimate in $L^2$? On torus, does estimate in $H^s$ with $s>0$ implies estimate in $L^2$?}
Let $s=0$ in (\ref{energy3}) and take sum in $q$, then we get
\begin{equation}\label{energy-l2}
\begin{split}
\|B(t)\|_{L^2}
\leq&\ \sum_{q\geq0}e^{-\mu\lambda_q^\alpha t} \|B_q(0)\|_{L^2}\\
&+ \sum_{q\geq0}\int_0^t e^{-\mu\lambda_q^\alpha (t-\tau)}\left\|\widetilde{\Delta_q}\left([B, \Delta_q]J_x\right)(\tau)\right\|_{L^2} \, d\tau\\
&+  \sum_{q\geq0}\int_0^t e^{-\mu\lambda_q^\alpha (t-\tau)}\left\| \widetilde{\Delta_q}\left([J, \Delta_q]B_x\right) (\tau)\right\|_{L^2}\, d\tau \\
&+ \sum_{q\geq0} \int_0^t e^{-\mu\lambda_q^\alpha (t-\tau)}\left\| \widetilde{\Delta_q}\left(BJ_{x,q}\right)(\tau) \right\|_{L^2}\, d\tau\\
& + \sum_{q\geq0}\int_0^t e^{-\mu\lambda_q^\alpha (t-\tau)} \left\| \widetilde{\Delta_q}\left(JB_{x,q}\right)(\tau) \right\|_{L^2} \, d\tau\\
=:&\ I_6+I_7+I_8+I_9+I_{10}.
\end{split}
\end{equation}
It is easy to see that $I_6\leq \|B(0)\|_{L^2}$. For $I_7$ we apply Lemma \ref{le-d1} with $\max\{2, \frac52-\alpha\}<m<\frac52$, $s_1=0$ and $s_2=-2$ to infer
\begin{equation}\notag
\begin{split}
I_7\leq & \sum_{q\geq 0} \int_0^t c_q\lambda_q^{-(m-\frac52)} e^{-\mu \lambda_q^\alpha (t-\tau)}\|B(\tau)\|_{L^2}\|B(\tau)\|_{H^m}\, d\tau\\
\lesssim &  \int_0^t (t-\tau)^{\frac{m-\frac52}{\alpha}}\|B(\tau)\|_{L^2}\|B(\tau)\|_{H^m}\, d\tau\\
\lesssim & \left(\sup_{t\in[0,T)}\|B(t)\|_{L^2} \right) \left(\sup_{t\in(0,T)}t^{\frac{m-(\frac52-\alpha)}{\alpha}}\|B(t)\|_{H^m} \right)\int_0^t (t-\tau)^{\frac{m-\frac52}{\alpha}}\tau^{\frac{\frac52-\alpha-m}{\alpha}}\, d\tau\\
\leq & C\left(\sup_{t\in[0,T)}\|B(t)\|_{L^2} \right) \left(\sup_{t\in(0,T)}t^{\frac{m-(\frac52-\alpha)}{\alpha}}\|B(t)\|_{H^m} \right)
\end{split}
\end{equation}
where we used the estimate of the time integral
\begin{equation}\notag
\int_0^t (t-\tau)^{\frac{m-\frac52}{\alpha}}\tau^{\frac{\frac52-\alpha-m}{\alpha}}\, d\tau
=\int_0^1 (1-\tau')^{\frac{m-\frac52}{\alpha}}(\tau')^{\frac{\frac52-\alpha-m}{\alpha}}\, d\tau'\leq C
\end{equation}
for the appropriate parameter $m$. Again in analogy with the estimate for $I_7$, we take $\max\{2, \frac52-\alpha\}<m<\frac52$ and apply Lemma \ref{le-d1} to $I_8$ with $s_1=s_2=-1$, apply Lemma \ref{le-d2} to $I_9$ with $s_1=0$ and $s_2=-2$, and apply Lemma \ref{le-d2} to $I_{10}$ with $s_1=s_2=-1$ to yield 
\begin{equation}\notag
I_8+I_9+I_{10}
\leq  C\left(\sup_{t\in[0,T)}\|B(t)\|_{L^2} \right) \left(\sup_{t\in(0,T)}t^{\frac{m-(\frac52-\alpha)}{\alpha}}\|B(t)\|_{H^m} \right).
\end{equation}
Therefore we have shown for $\max\{2, \frac52-\alpha\}<m<\frac52$
\begin{equation}\notag
\sup_{t\in[0,T)}\|B(t)\|_{L^2}
\leq  \|B_0\|_{L^2}+C\left(\sup_{t\in[0,T)}\|B(t)\|_{L^2} \right) \left(\sup_{t\in(0,T)}t^{\frac{m-(\frac52-\alpha)}{\alpha}}\|B(t)\|_{H^m} \right).
\end{equation}
The lower bound 
\begin{equation}\notag
\begin{split}
&\sup_{t\in[0,T)}\|B(t)\|_{L^2}\\
\geq&\  \left(\sum_{q\geq0}e^{-2\mu\lambda_q^\alpha t} \|B_q(0)\|^2_{L^2}\right)^{\frac12}\\
&-C\left(\sup_{t\in[0,T)}\|B(t)\|_{L^2} \right) \left(\sup_{t\in(0,T)}t^{\frac{m-(\frac52-\alpha)}{\alpha}}\|B(t)\|_{H^m} \right)
\end{split}
\end{equation}
can be established in an analogous way. We omit the details here. 

\medskip

\subsection{Proof of Theorem \ref{thm-local} and Theorem \ref{thm-local1}}
\label{sec-proof1}

We consider the approximating system 
\begin{equation}\label{sys-app}
\begin{split}
B^{k}_t+aB^{k-1} J_x^{k}+bJ^{k-1}B_x^{k}+\mu \Lambda^\alpha B^{k}=&\ 0,\\
B_x^k=&\ \mathcal H J^k
\end{split}
\end{equation}
for $k\geq 0$, and $B^{-1}=J^{-1}\equiv 0$, with initial data $B^{k}(x,0)=B_0(x)$ for all $k\geq 0$. In analogy with establishing estimates (\ref{energy4})-(\ref{energy6}) we obtain for $\frac52-\alpha<s<\frac12$ and $k\geq 0$,
\begin{equation}\label{energy-app1}
\begin{split}
&\|B^k\|_{L^{\frac{\alpha}{(s+\alpha-\frac52)}}(0,T; H^s)}\\
\leq&\ C(T) \|B(0)\|_{H^{\frac{5}{2} -\alpha}}
+C \|B^{k-1}\|_{L^{\frac{\alpha}{(s+\alpha-\frac52)}}(0,T; H^s)} \|B^{k}\|_{L^{\frac{\alpha}{(s+\alpha-\frac52)}}(0,T; H^s)}
\end{split}
\end{equation} 
\begin{equation}\label{energy-app2}
\begin{split}
&\|B^k\|_{L^{2}(0,T; H^{\frac52-\frac{\alpha}{2}})}\\
\leq&\ C(T) \|B(0)\|_{H^{\frac{5}{2} -\alpha}}+C \|B^{k-1}\|_{L^{\frac{\alpha}{(s+\alpha-\frac52)}}(0,T; H^s)}\|B^{k}\|_{L^{\frac{\alpha}{(s+\alpha-\frac52)}}(0,T; H^s)}
\end{split}
\end{equation}
\begin{equation}\label{energy-app3}
\begin{split}
&\|B^k\|_{L^{\infty}(0,T; H^{\frac52-\alpha})}\\
\leq &\ \|B(0)\|_{H^{\frac{5}{2} -\alpha}}+C \|B^{k-1}\|_{L^{2}(0,T; H^{\frac52-\frac{\alpha}{2}})}\|B^{k}\|_{L^{2}(0,T; H^{\frac52-\frac{\alpha}{2}})}.
\end{split}
\end{equation}
The estimate (\ref{energy-app1}) indicates that for $T>0$ sufficiently small
\begin{equation}\label{energy-app4}
\|B^k\|_{L^{\frac{\alpha}{(s+\alpha-\frac52)}}(0,T; H^s)}
\leq 2C(T) \|B(0)\|_{H^{\frac{5}{2} -\alpha}} \ \ \ \forall \ \ k\geq0.
\end{equation} 
It then follows from (\ref{energy-app2}) and (\ref{energy-app3}) that for even smaller $T>0$ 
\begin{equation}\label{energy-app5}
\begin{split}
\|B^k\|_{L^2(0,T; H^{\frac52-\frac{\alpha}2})}\leq&\ 2C(T) \|B(0)\|_{H^{\frac{5}{2} -\alpha}} \\
\|B^k\|_{L^\infty(0,T; H^{\frac52-\alpha})}\leq&\ (1+C(T)) \|B(0)\|_{H^{\frac{5}{2} -\alpha}} 
\end{split}
\end{equation} 
for all $k\geq 0$. In view of (\ref{energy-low}) and the upper bound in (\ref{energy-app5}), it follows from the dominated convergence theorem that
\begin{equation}\label{energy-app6}
\inf_{t\in[0,T)} \|B^k(t)\|_{H^{\frac52-\alpha}}\geq (1-C(T)) \|B_0\|_{H^{\frac52-\alpha}} \ \ \ \mbox{uniformly in} \ \ k.
\end{equation}
From the proof of (\ref{est-b1}) we also see for any $\beta>0$ and $k\geq 0$
\begin{equation}\label{energy-app7}
 \|t^{\frac{\beta}{\alpha}}B^k(t)\|_{L^\infty(0,T; H^{\frac52-\alpha+\beta})}\leq C(T, \beta) \|B_0\|_{H^{\frac52-\alpha}}.
\end{equation}
By the analysis of Subsection \ref{sec-l2} we also have the uniform estimates
\begin{equation}\label{energy-app8}
\begin{split}
\sup_{t\in[0,T)} \|B^k(t)\|_{L^2} \leq (1+C(T)) \|B_0\|_{L^2},\\
\inf_{t\in[0,T)} \|B^k(t)\|_{L^2} \leq (1-C(T)) \|B_0\|_{L^2}.
\end{split}
\end{equation}

The next step is to show that $\{B^k\}$ is a Cauchy sequence in the Banach space $L^\infty(0,T; L^2)$. We consider $\tilde B^k=B^k-B^{k-1}$ which satisfies 
\begin{equation}\label{eq-diff}
\tilde B^{k+1}_t+B^k \tilde J^{k+1}_x+\tilde B^k J^k_x-J^k\tilde B^{k+1}_x-\tilde J^k B_x^k+\mu \Lambda^\alpha \tilde B^{k+1}=0
\end{equation}
with $B^k_x=\mathcal H J^k$ and $\tilde J^k=J^{k}-J^{k-1}$ for $k\geq 0$ and $B^{-1}\equiv J^{-1} \equiv 0$. Analogous arguments as in Subsection \ref{sec-priori} lead to
\begin{equation}\notag
\begin{split}
\|\tilde B^{k+1}_q(t)\|_{L^2}\leq& \int_0^te^{-\mu\lambda_q^\alpha (t-\tau)} \| \widetilde \Delta_q [B^k, \Delta_q] \tilde J_x^{k+1}(\tau)\|_{L^2}\, d\tau\\
&+ \int_0^te^{-\mu\lambda_q^\alpha (t-\tau)} \| \widetilde \Delta_q \left(B^k\tilde J_{x,q}^{k+1}\right)(\tau)\|_{L^2}\, d\tau\\
&-\int_0^te^{-\mu\lambda_q^\alpha (t-\tau)} \| \widetilde \Delta_q [J^k, \Delta_q] \tilde B_x^{k+1}(\tau)\|_{L^2}\, d\tau\\
&- \int_0^te^{-\mu\lambda_q^\alpha (t-\tau)} \| \widetilde \Delta_q \left(J^k\tilde B_{x,q}^{k+1}\right)(\tau)\|_{L^2}\, d\tau\\
&- \int_0^te^{-\mu\lambda_q^\alpha (t-\tau)} \|\Delta_q \left(\tilde B^k J_{x}^{k}\right)(\tau)\|_{L^2}\, d\tau\\
&+ \int_0^te^{-\mu\lambda_q^\alpha (t-\tau)} \|\Delta_q \left(\tilde J^k B_{x}^{k}\right)(\tau)\|_{L^2}\, d\tau.
\end{split}
\end{equation}
In the order it appears, we apply Lemmas \ref{le-d1}, \ref{le-d2} and \ref{le-d3} with appropriate norms to the integrals on the right hand side of the inequality above\\
(i) Lemma \ref{le-d1} with $m=0$, $s_1=\frac52-\alpha+\beta<\frac32$ and $s_2=-2$,\\
(ii) Lemma \ref{le-d2} with $m=0$, $s_1=\frac52-\alpha+\beta<\frac12$ and $s_2=-2$,\\
(iii) Lemma \ref{le-d1} with $m=0$, $s_1=\frac32-\alpha+\beta<\frac32$ and $s_2=-1$,\\
(iv) Lemma \ref{le-d2} with $m=0$, $s_1=\frac32-\alpha+\beta<\frac12$ and $s_2=-1$,\\
(v) Lemma \ref{le-d3} with $s_1=0$ and $s_2=\frac12-\alpha+\beta<\frac12$,\\
(vi) Lemma \ref{le-d3} with $s_1=-1$ and $s_2=\frac32-\alpha+\beta<\frac12$,\\
and we obtain
\begin{equation}\notag
\begin{split}
&\sum_{q\geq -1}\|\tilde B^{k+1}_q(t)\|_{L^2}\\
\lesssim& \sum_{q\geq -1}\int_0^tc_q \lambda_q^{\alpha-\beta}e^{-\mu\lambda_q^\alpha (t-\tau)}\|B^k(\tau)\|_{H^{\frac52-\alpha+\beta}}\left(\|\tilde B^k(\tau)\|_{L^2}+\|\tilde B^{k+1}(\tau)\|_{L^2} \right)\, d\tau\\
\lesssim & \int_0^t(t-\tau)^{-1+\frac{\beta}{\alpha}}\|B^k(\tau)\|_{H^{\frac52-\alpha+\beta}}\left(\|\tilde B^k(\tau)\|_{L^2}+\|\tilde B^{k+1}(\tau)\|_{L^2} \right)\, d\tau\\
\leq &\ C_0 \|t^{\frac{\beta}{\alpha}} B^k\|_{L^\infty(0,T; H^{\frac52-\alpha+\beta})}\left(\|\tilde B^k\|_{L^\infty(0,T; L^2)}+\|\tilde B^{k+1}\|_{L^\infty(0,T; L^2)} \right)\\
&\cdot \int_0^t(t-\tau)^{-1+\frac{\beta}{\alpha}}\tau^{-\frac{\beta}{\alpha}}\, d\tau
\end{split}
\end{equation}
for a constant $C_0>0$.
Note that the time integral is bounded if $0<\frac{\beta}{\alpha}<1$ which is satisfied upon conditions above to apply the Lemmas. Hence we have for some $0<\beta<\alpha$
\begin{equation}\notag
\|\tilde B^{k+1}\|_{L^\infty(0,T; L^2)}\leq C_0 \|t^{\frac{\beta}{\alpha}} B^k\|_{L^\infty(0,T; H^{\frac52-\alpha+\beta})}\left(\|\tilde B^k\|_{L^\infty(0,T; L^2)}+\|\tilde B^{k+1}\|_{L^\infty(0,T; L^2)} \right).
\end{equation}
Again for sufficiently small $T>0$ with 
\[\|t^{\frac{\beta}{\alpha}} B^k\|_{L^\infty(0,T; H^{\frac52-\alpha+\beta})}\leq \frac{1}{4C_0},\]
we deduce 
\begin{equation}\notag
\|\tilde B^{k+1}\|_{L^\infty(0,T; L^2)}\leq \frac13\|\tilde B^k\|_{L^\infty(0,T; L^2)}
\end{equation}
which implies by iteration
\begin{equation}\label{energy-app9}
\|\tilde B^{k+1}\|_{L^\infty(0,T; L^2)}\leq 3^{-k} C.
\end{equation}
Hence $\{B^k\}$ is a Cauchy sequence in $L^\infty(0,T; L^2)$ which converges to $B\in L^\infty(0,T; L^2)$.  It follows from 
(\ref{energy-app5}) and (\ref{energy-app7}) that
\begin{equation}\label{energy-app10}
\|t^{\frac{\beta}{\alpha}} B\|_{L^\infty(0,T; H^{\frac52-\alpha+\beta})}\leq C, \ \ \ \beta\geq0,
\end{equation}
\begin{equation}\label{energy-app11}
\lim_{t\to 0}\|t^{\frac{\beta}{\alpha}} B\|_{H^{\frac52-\alpha+\beta}}=0, \ \ \ \beta>0.
\end{equation}
We claim that the convergence of $\{B^k\}$ to $B$ also occurs in Sobolev spaces with higher regularity. Indeed, interpolation of (\ref{energy-app7}) and (\ref{energy-app9}) gives 
\begin{equation}\notag
 \|t^{\frac{s\beta}{\alpha}}B^k(t)\|_{L^\infty(0,T; H^{s(\frac52-\alpha+\beta)})}\leq c(T, \beta) 3^{-k(1-s)}
\end{equation}
for $s\in[0,1)$ and $\beta\geq 0$. It is thus obvious that $\{t^{\frac{\beta'}{\alpha'}}B^k\}$ is a Cauchy sequence in $L^\infty(0,T; H^{\beta'})$ for any $\beta'\geq 0$ and $\alpha'>\alpha$, and the sequence converges to $B$ in such spaces. It then follows from Sobolev embedding theorem that the sequence $\{B^k\}$ is smooth and the derivatives of the sequence converge to the derivatives of $B$ for $t\in(0,T)$ uniformly in $x$. Therefore (\ref{energy-app8}) implies 
\[B\in L^\infty([0, T); H^{\frac52-\alpha})\cap L^2([0,T); H^{\frac52-\frac{\alpha}{2}})\]
and $B$ is a classical solution of (\ref{emhd-1d}).

In the end we show the continuity in time of $B$, that is
\[B\in C([0,T); H^{\frac52-\alpha}).\]
The uniform upper and lower bounds (\ref{energy-app5}), (\ref{energy-app6}) and (\ref{energy-app8}) imply the continuity at $t=0^+$. On the other hand, for any $\delta\in (0,T)$ it follows from the estimate (\ref{energy-app10}) that
\[B_t=-BJ_x+B_xJ-\mu\Lambda^\alpha B\in L^1([\delta, T); H^{\frac52-\alpha})\]
which implies 
\[B\in C([\delta,T); H^{\frac52-\alpha}).\]

We point out that slight modification of the analysis for the difference equation (\ref{eq-diff}) gives the uniqueness of the solution. Thus the proof of Theorem \ref{thm-local} is complete. 

Note that when $a=0$ and $b=1$, the requirement on the the parameters for the a priori estimates is $\frac52-\alpha<s<\frac32$, which implies $\alpha>1$. Hence we have also proved Theorem \ref{thm-local1}. 

\cbdu

\bigskip

\section{Local well-posedness in the subcritical regime}
\label{sec-local2}

This section is devoted to a proof of Theorem \ref{thm-local2}, which is analogous to the analysis of Section \ref{sec-local}.
The major task is to establish the a priori estimates.  
Similar treatment as before, we have from (\ref{energy1}) for any $0\leq s<\frac12$
\begin{equation}\label{energy1-lw}
\begin{split}
\lambda_q^{s}\|B_q(t)\|_{L^2}
\leq&\ \lambda_q^{s}e^{-\mu\lambda_q^\alpha t} \|B_q(0)\|_{L^2}\\
&+ \int_0^t \lambda_q^{s}e^{-\mu\lambda_q^\alpha (t-\tau)}\left\|\widetilde{\Delta_q}\left([B, \Delta_q]J_x\right)(\tau)\right\|_{L^2} \, d\tau\\
&+  \int_0^t \lambda_q^{s}e^{-\mu\lambda_q^\alpha (t-\tau)}\left\| \widetilde{\Delta_q}\left([J, \Delta_q]B_x\right) (\tau)\right\|_{L^2}\, d\tau \\
&+  \int_0^t \lambda_q^{s}e^{-\mu\lambda_q^\alpha (t-\tau)}\left\| \widetilde{\Delta_q}\left(BJ_{x,q}\right)(\tau) \right\|_{L^2}\, d\tau\\
& + \int_0^t \lambda_q^{s}e^{-\mu\lambda_q^\alpha (t-\tau)} \left\| \widetilde{\Delta_q}\left(JB_{x,q}\right)(\tau) \right\|_{L^2} \, d\tau.
\end{split}
\end{equation}
Lemma \ref{le-m1} and Lemma \ref{le-m2} together imply that 
\begin{equation}\notag%\label{term1-lw}
\begin{split}
&\left\| \lambda_q^{s}e^{-\mu\lambda_q^\alpha t} \|B_q(0)\|_{L^{2}}\right\|_{L^{\frac{\alpha}{s}}(0,T; l^2)}\\
\leq &\ \left\| e^{-\mu\lambda_q^\alpha t} B_q(0)\right\|_{L^{\frac{\alpha}{s}}(0,T; H^{s})}\\
\leq&\ C(T) \|B(0)\|_{L^{2}}
\end{split}
\end{equation}
where the constant $C(T)\to 0$ as $T\to 0$.

We apply Lemma \ref{le-d1} with $m=0$, $s_1=s<\frac32$ and $s_2=s-2<\frac12$ to the second term on the right hand side of (\ref{energy1-lw}) 
\begin{equation}\notag
\begin{split}
& \int_0^t \lambda_q^{s}e^{-\mu\lambda_q^\alpha (t-\tau)}\left\|\widetilde{\Delta_q}\left([B, \Delta_q]J_x\right)(\tau)\right\|_{L^2} \, d\tau\\
\leq & \int_0^t c_q\lambda_q^{\frac52-s}e^{-\mu\lambda_q^\alpha (t-\tau)} \|B\|_{H^{s}}\|J_x\|_{H^{s-2}} \, d\tau\\
%\leq & \int_0^t c_q\lambda_q^{\frac{5}{2}-s}e^{-\mu\lambda_q^\alpha (t-\tau)} \|B\|^2_{H^{s}} \, d\tau\\
\leq & \int_0^t c_q(t-\tau)^{-\frac{1}{\alpha}(\frac52-s)} \|B\|_{H^{s}}^2 \, d\tau\\
\leq &\ t^{1-\frac{5}{2\alpha}}\int_0^t c_q(t-\tau)^{-1+\frac{s}{\alpha}} \|B\|_{H^{s}}^2 \, d\tau\\
\end{split}
\end{equation}
%where we have used the fact $x^ae^{-x}\leq C$ for $a\geq 0$ and 
for $\alpha>\frac52$. Thus we apply Hardy-Littlewood-Sobolev inequality again
\begin{equation}\label{term1-lw}
\begin{split}
&\left\| \int_0^t \lambda_q^{s}e^{-\mu\lambda_q^\alpha (t-\tau)}\left\|\widetilde{\Delta_q}\left([B, \Delta_q]J_x\right)(\tau)\right\|_{L^2} \, d\tau\right\|_{L^{\frac{\alpha}{s}}(0,T; l^2)}\\
\lesssim &\ T^{1-\frac{5}{2\alpha}} \left\| \int_0^t (t-\tau)^{-1+\frac{s}{\alpha}} \|B\|_{H^{s}}^2 \, d\tau\right\|_{L^{\frac{\alpha}{s}}(0,T)}\\
\lesssim &\ T^{1-\frac{5}{2\alpha}}\left\|\|B\|_{H^{s}}^2\right\|_{L^{\frac{\alpha}{2s}}(0,T)}\\
\lesssim &\ T^{1-\frac{5}{2\alpha}}\|B\|^2_{L^{\frac{\alpha}{s}}(0,T; H^{s})}.
\end{split}
\end{equation}

Analogously, we have
\begin{equation}\label{energy2-lw}
\begin{split}
\lambda_q^{\frac{\alpha}{2}}\|B_q(t)\|_{L^2}
\leq&\ \lambda_q^{\frac{\alpha}{2}}e^{-\mu\lambda_q^\alpha t} \|B_q(0)\|_{L^2}\\
&+ \int_0^t \lambda_q^{\frac{\alpha}{2}}e^{-\mu\lambda_q^\alpha (t-\tau)}\left\|\widetilde{\Delta_q}\left([B, \Delta_q]J_x\right)(\tau)\right\|_{L^2} \, d\tau\\
&+  \int_0^t \lambda_q^{\frac{\alpha}{2}}e^{-\mu\lambda_q^\alpha (t-\tau)}\left\| \widetilde{\Delta_q}\left([J, \Delta_q]B_x\right) (\tau)\right\|_{L^2}\, d\tau \\
&+  \int_0^t \lambda_q^{\frac{\alpha}{2}}e^{-\mu\lambda_q^\alpha (t-\tau)}\left\| \widetilde{\Delta_q}\left(BJ_{x,q}\right)(\tau) \right\|_{L^2}\, d\tau\\
& + \int_0^t \lambda_q^{\frac{\alpha}{2}}e^{-\mu\lambda_q^\alpha (t-\tau)} \left\| \widetilde{\Delta_q}\left(JB_{x,q}\right)(\tau) \right\|_{L^2} \, d\tau.
\end{split}
\end{equation}
Again in view of Lemma \ref{le-m1} and Lemma \ref{le-m2}, we obtain
\begin{equation}\notag
\begin{split}
&\left\| \lambda_q^{\frac{\alpha}{2}}e^{-\mu\lambda_q^\alpha t} \|B_q(0)\|_{L^{2}}\right\|_{L^{2}(0,T; l^2)}\\
\leq &\ \left\| e^{-\mu\lambda_q^\alpha t} B_q(0)\right\|_{L^{2}(0,T; H^{\frac{\alpha}{2}})}\\
\leq&\ C(T) \|B(0)\|_{L^{2}}
\end{split}
\end{equation}
where the constant $C(T)\to 0$ as $T\to 0$.

We apply Lemma \ref{le-d1} with $m=0$, $s_1=s<\frac32$ and $s_2=s-2<\frac12$ to the second term of the right hand side of (\ref{energy2-lw}) 
\begin{equation}\notag
\begin{split}
& \int_0^t \lambda_q^{\frac{\alpha}{2}}e^{-\mu\lambda_q^\alpha (t-\tau)}\left\|\widetilde{\Delta_q}\left([B, \Delta_q]J_x\right)(\tau)\right\|_{L^2} \, d\tau\\
\leq & \int_0^t c_q\lambda_q^{\frac{\alpha}{2}+\frac52-2s}e^{-\mu\lambda_q^\alpha (t-\tau)} \|B\|_{H^{s}}\|J_x\|_{H^{s-2}} \, d\tau\\
\leq & \int_0^t c_q(t-\tau)^{-\frac{1}{\alpha}(\frac{\alpha}{2}+\frac52-2s)} \|B\|_{H^{s}}^2 \, d\tau\\
\leq &\ t^{1-\frac{5}{2\alpha}}\int_0^t c_q(t-\tau)^{-\frac{3}{2}+\frac{2s}{\alpha}} \|B\|_{H^{s}}^2 \, d\tau.
\end{split}
\end{equation}
%where we have used the fact $x^ae^{-x}\leq C$ for $a\geq 0$ and $\alpha>\frac52$. 
It follows from Hardy-Littlewood-Sobolev inequality that
\begin{equation}\label{term2-lw}
\begin{split}
&\left\| \int_0^t \lambda_q^{\frac{\alpha}{2}}e^{-\mu\lambda_q^\alpha (t-\tau)}\left\|\widetilde{\Delta_q}\left([B, \Delta_q]J_x\right)(\tau)\right\|_{L^2} \, d\tau\right\|_{L^{2}(0,T; l^2)}\\
\lesssim &\ T^{1-\frac{5}{2\alpha}} \left\| \int_0^t (t-\tau)^{-\frac32+\frac{2s}{\alpha}} \|B\|_{H^{s}}^2 \, d\tau\right\|_{L^{2}(0,T)}\\
\lesssim &\ T^{1-\frac{5}{2\alpha}}\left\|\|B\|_{H^{s}}^2\right\|_{L^{\frac{\alpha}{2s}}(0,T)}\\
\lesssim &\ T^{1-\frac{5}{2\alpha}}\|B\|^2_{L^{\frac{\alpha}{s}}(0,T; H^{s})}.
\end{split}
\end{equation}

The other three integrals in (\ref{energy1-lw}) and (\ref{energy3-lw}) can be handled in a similar way. We point out that the estimate for the forth term on the right hand side of (\ref{energy1-lw}) and (\ref{energy2-lw}) requires $s<\frac12$.
Thus in view of (\ref{term1-lw}) and (\ref{term2-lw}) we claim for $0\leq s<\frac12$
\begin{equation}\label{energy3-lw}
\|B\|_{L^{\frac{\alpha}{s}}(0,T; H^s)}\leq C(T) \|B(0)\|_{L^{2}}+C T^{1-\frac{5}{2\alpha}}\|B\|^2_{L^{\frac{\alpha}{s}}(0,T; H^{s})}
\end{equation}
\begin{equation}\label{energy4-lw}
\|B\|_{L^{2}(0,T; H^{\frac{\alpha}{2}})}\leq C(T) \|B(0)\|_{L^{2}}+C T^{1-\frac{5}{2\alpha}}\|B\|^2_{L^{\frac{\alpha}{s}}(0,T; H^{s})}.
\end{equation}
Similarly as before, we can show 
\begin{equation}\label{energy5-lw}
\|B\|_{L^{\infty}(0,T; L^{2})}\leq C(T) \|B(0)\|_{L^{2}}+C T^{1-\frac{5}{2\alpha}}\|B\|^2_{L^{\frac{\alpha}{s}}(0,T; H^{s})},
\end{equation}
\begin{equation}\label{energy6-lw}
\|B(t)\|_{L^{2}}\geq \left(\sum_{q\geq -1}e^{-2\mu\lambda_q^\alpha t} \|B_q(0)\|^2_{L^{2}}\right)^{\frac12}-C t^{1-\frac{5}{2\alpha}}\|B\|^2_{L^{\frac{\alpha}{s}}(0,t; H^{s})}.
\end{equation}
Following the lines of Subsection \ref{sec-smooth}, we can obtain the smoothing estimates (\ref{est-b3}) and (\ref{est-b4}).
With the estimates (\ref{energy3-lw})-(\ref{energy6-lw}) and (\ref{est-b3})-(\ref{est-b4}) at hand, applying a similar analysis as in Subsection \ref{sec-proof1} shows that there is a unique solution $B(t)$ of (\ref{emhd-1d}) with $\alpha>\frac52$ in $C\left([0,T);   L^2(\mathbb S^1)\right)\cap L^2\left([0,T);  H^{\frac{\alpha}2}(\mathbb S^1)\right)$.

\bigskip

\section{Analytic solutions for small $\alpha$}
In the previous sections we reveal obstacles to obtain local well-posedness for the equations (\ref{emhd-1d}) and (\ref{e1d3}) with small $\alpha$ in Sobolev spaces. The obstacles are related to the loss of derivative in the models, which comes from the very singular structure of the nonlinear terms. Nevertheless, we are able to show the existence of local analytic solutions to (\ref{e1d3}) which is the main purpose of this section. In order to prove Theorem \ref{thm-analytic} we only show the a priori estimates which establishes the analyticity; the existence follows from standard approximating arguments. 

Denote the projection $B^N(x,t)=\mathbb P_{\leq N} B(x,t)$ and $J^N=-\Lambda B^N$ for any integer $N\geq 1$. Acting the projection $\mathbb P_{\leq N}$ on equation (\ref{e1d3}) and taking Fourier transform yields
\begin{equation}\notag
\widehat B_t^N(k,t)=i\sum_{m+n=k, |m|, |n|,|k|\leq N} m|n|\widehat B(m) \widehat B(n)-\mu |k|^\alpha\widehat B(k).
\end{equation}
Denote 
\[\xi^N(k,t)=\widehat B^N(k,t) e^{\frac12\mu|k|^\alpha t}, \ \ Y^N(t)=\sum_{|k|\leq N} |k|^8 |\xi^N(k,t)|^2.\]
The goal is to show that $Y^N(t)$ satisfies a Riccati inequality. Note that $\bar \xi^N(k)=\xi^N(-k)$ since $B(x,t)$ is real-valued. First we have
\begin{equation}\label{est1-ana}
\begin{split}
\frac{d}{dt} \xi^N(k,t)=&\ i \sum_{\substack{m+n=k, \\ |m|, |n|,|k|\leq N}} m|n| e^{\frac12\mu(|k|^\alpha-|m|^\alpha-|n|^\alpha)t}\xi^N(m) \xi^NB(n)\\
&-\frac12\mu |k|^\alpha\xi^N(k).
\end{split}
\end{equation}
Denote $\gamma:=\gamma_{m,n,k}=\frac12\mu(|m|^\alpha+|n|^\alpha-|k|^\alpha)$. It then follows from (\ref{est1-ana}) that
\begin{equation}\label{est2-ana}
\begin{split}
\frac{d}{dt} Y^N(t)=&\sum_{|k|\leq N}|k|^8 \xi^N(k)\frac{d}{dt}\bar \xi^N(k)+\sum_{|k|\leq N}|k|^8 \bar\xi^N(k)\frac{d}{dt}\xi^N(k)\\
=&\ 2\Re\left(i\sum_{\substack{m+n+k=0, \\ |m|, |n|,|k|\leq N}} e^{-\gamma t} m|n||k|^8\xi^N(m)\xi^N(n)\xi^N(k)\right)\\
&-\mu |k|^{\alpha+8}|\xi^N(k)|^2\\
=&\ 2\Re\left(i\sum_{\substack{m+n+k=0, \\ |m|, |n|,|k|\leq N}} m|n||k|^8\xi^N(m)\xi^N(n)\xi^N(k)\right)\\
&+ 2\Re\left(i\sum_{\substack{m+n+k=0, \\ |m|, |n|,|k|\leq N}} (e^{-\gamma t}-1)m|n||k|^8\xi^N(m)\xi^N(n)\xi^N(k)\right)\\
&-\mu \sum_{|k|\leq N}|k|^{\alpha+8}|\xi^N(k)|^2\\
=:&\ I_1+I_2+I_3.
\end{split}
\end{equation}
In view of symmetrization we can write
%and the relation $m+n+k=0$ we infer
\begin{equation}\notag
I_1=\frac23 \Re\left( i\sum_{\substack{m+n+k=0, \\ |m|, |n|,|k|\leq N}}(m|n||k|^8+m|k||n|^8+k|n||m|^8)\xi^N(m)\xi^N(n)\xi^N(k)\right).
\end{equation}
Hence we infer by using the relation $m+n+k=0$, Young's convolution inequality and H\"older's inequality 
\begin{equation}\notag
\begin{split}
|I_1|\lesssim & \sum_{\substack{m+n+k=0, \\ |m|\leq |n|\leq |k|\leq N}}\left|m|n||k|^8+m|k||n|^8+k|n||m|^8\right| \left|\xi^N(m)\right|\left|\xi^N(n)\right|\left|\xi^N(k)\right|\\
\lesssim & \sum_{\substack{m+n+k=0, \\ |m|\leq |n|\leq |k|\leq N}}|m|^2|n|^4|k|^4 \left|\xi^N(m)\right|\left|\xi^N(n)\right|\left|\xi^N(k)\right|\\
\lesssim &\ Y^N\sum_{|m|\leq N}|m|^2 \left|\xi^N(m)\right|\\
\lesssim &\ (Y^N)^{\frac32}.
\end{split}
\end{equation}

Applying the facts $1-e^{-\gamma t}\leq \gamma t$ for $\gamma>0$ and $0<\gamma<\min\{|m|, |n|\}$, the relation $m+n+k=0$ again, Young's and H\"older's inequalities we have
\begin{equation}\notag
\begin{split}
|I_2|\leq &\ 2\sum_{\substack{m+n+k=0, \\ |m|, |n|,|k|\leq N}} \gamma t |m||n||k|^8\left|\xi^N(m)\right|\left|\xi^N(n)\right|\left|\xi^N(k)\right|\\
\lesssim &\ \mu t\sum_{\substack{m+n+k=0, \\ |m|\leq |n|\leq |k|\leq N}}  |m|^2|n||k|^8\left|\xi^N(m)\right|\left|\xi^N(n)\right|\left|\xi^N(k)\right|\\
\lesssim &\ \mu t\sum_{\substack{m+n+k=0, \\ |m|\leq |n|\leq |k|\leq N}}  |m|^{3-\alpha}|n|^{4+\frac{\alpha}{2}}|k|^{4+\frac{\alpha}{2}}\left|\xi^N(m)\right|\left|\xi^N(n)\right|\left|\xi^N(k)\right|\\
\lesssim &\ \mu t\left(\sum_{|m|\leq N}|m|^{3-\alpha} \left|\xi^N(m)\right|\right)\left(\sum_{|k|\leq N}|k|^{8+\alpha} \left|\xi^N(k)\right|^2\right)\\
\lesssim &\ \mu t\left(\sum_{|m|\leq N}|m|^{8} \left|\xi^N(m)\right|^2\right)^{\frac12}\left(\sum_{|k|\leq N}|k|^{8+\alpha} \left|\xi^N(k)\right|^2\right).
\end{split}
\end{equation}
Inserting the estimates of $I_1$ and $I_2$ into (\ref{est2-ana}) we obtain
\begin{equation}\label{est3-ana}
\frac{d}{dt} Y^N(t)\leq C_1 (Y^N)^{\frac32}+C_2\mu\left((Y^N)^{\frac12}t-1\right) \left(\sum_{|k|\leq N}|k|^{8+\alpha} \left|\xi^N(k)\right|^2\right)
\end{equation}
for some absolute constants $C_1, C_2>0$. Note that $Y^N(0)=\|B_0^N\|_{L^2}^2$. It follows from (\ref{est3-ana}) that $Y^N$ is bounded from above on $[0,T)$ for a small time $T>0$, and the bound is uniform in $N$. It implies the limit function $B$ of the sequence $\{B^N\}$ is analytic on $[0,T)$. 

\bigskip

\section{Singularity formation}

This section provides a proof for Theorem \ref{thm-s1}. % and \ref{thm-s2}.
We first show that the maximum principle holds for (\ref{e1d3}) with $\mu\geq 0$ and $0< \alpha<2$. 
\begin{Lemma}\label{le-max}
Let $B_0\in L^\infty$ be the initial data of (\ref{e1d3}) with $0< \alpha<2$. The solution $B(t)$ satisfies
\[ \|B(t)\|_{L^\infty}\leq \|B_0\|_{L^\infty} \ \ \ \forall \ \ t\geq 0.\]
\end{Lemma}
\pf
Formally, let $B(t)$ be a smooth solution of (\ref{e1d3}) on $[0,T)$ with the initial data $B_0$. Let $0<t_1<t_2<T$ and $\varepsilon>0$. Denote $Q=[t_1, t_2]\times \mathbb S^1$. Note $B\in L^2(Q)$. We consider $\bar B=B-\varepsilon t$. Note that $\bar B$ is also smooth on $Q$ and its derivatives are bounded. There exists a point $(x^*, t^*)\in Q$ such that 
\[\bar B(x^*, t^*) =\sup_{(x,t)\in Q} \bar B(x,t).\]
Recall for $0<\alpha<2$, 
\[\Lambda^{\alpha} B(x,t)=\frac{1}{\pi} P.V.\int_{\mathbb R} \frac{B(x,t)-B(y,t)}{|x-y|^{1+\alpha}}\, dy\]
or
\[\Lambda^{\alpha} B(x,t)=\frac{1}{2\pi} P.V.\int_{\mathbb S^1} \frac{B(x,t)-B(y,t)}{|\sin \frac{x-y}{2}|^{1+\alpha}}\, dy,\]
%({\color{blue} Check out whether the formula is true for arbitrary $\alpha$, refer to Cordoba-Cordoba \cite{CC}. $\Lambda^2 B=-B_{xx}$; when $\alpha>2$, $\Lambda^\alpha B=-\Lambda^{\alpha-2} B_{xx}$. It might be true only for $0<\alpha<1$. })
and hence 
\[\Lambda^{\alpha} B(x^*,t^*)\geq 0.\]
At the maximum point $(x^*, t^*)$, we have 
\[\bar B_t(x^*, t^*)=B_t(x^*, t^*)-\varepsilon=B_x(x^*, t^*) J(x^*, t^*)-\mu \Lambda^\alpha B_t(x^*, t^*)-\varepsilon<0. \]
Since $(x^*, t^*)$ is the maximum point on $Q$, $t^*$ must be $t_1$. This argument is valid for arbitrary $\varepsilon>0$. Therefore,
\[\sup_{x\in \mathbb S^1} B(x,t_1)= \sup_{(x,t)\in Q} B(x,t).\]
An analogous argument also gives 
\[\inf_{x\in \mathbb S^1} B(x,t_1)= \inf_{(x,t)\in Q} B(x,t).\]
Since one can choose $t_1$ and $t_2$ arbitrarily, it completes the proof of the lemma. 

\cbdu

\begin{Remark}
As a special case of Lemma \ref{le-max}, solution starting from non-negative (positive) initial data remains non-negative (positive). Consequently, we can show
\[\|B(t)\|_{L^1}= \|B_0\|_{L^1}, \ \ \forall \ \ t\geq 0.\] 
Indeed, since $B(t)\geq 0$, we have
\begin{equation}\notag
\frac{d}{dt}\int |B|\, dx=\frac{d}{dt}\int B\, dx=\int B_xJ\, dx-\mu\int \Lambda^\alpha B\, dx 
\end{equation}
and note
\[\int B_xJ\, dx=\int J\mathcal H J\, dx=0, \]
by applying integration by parts, 
\[\int \Lambda^\alpha B\, dx=0.\]
Therefore, 
\begin{equation}\notag
\frac{d}{dt}\int |B|\, dx=0.
\end{equation}
Combining the maximum principle and the bound of $L^1$ norm, we also have the $L^2$ norm bounded.  
\end{Remark}

\medskip

{\textbf{Proof of Theorem \ref{thm-s1}:}}
Since $0\leq B_0\leq 1$ on $[0,\infty)$, it follows from Lemma \ref{le-max} that the solution satisfies $0\leq B(x, t)\leq 1$ for all $t\geq0$. Thanks to the properties of the Hilbert transform we know $B(\cdot,t)$ remains odd in $x$. Hence $B(0, t)=0$ for all $t\geq 0$, and $B(x,t)$ is compactly supported on $[-1,1]$ for all the time. 
%``Under the assumptions on the initial data $J_0$, we know $J(t,x)$ will remain positive (given the transport character of equation) and symmetric. Then $\mathcal H J$ will be antisymmetric and positive for $x\geq L$. This implies the following properties for $J(t,x)$: \\
%(i) $\sup J(t,x)\subset [-L, L]$;\\
%(ii) $\max_x J=J(0,0)=1$;\\
%(iii) $\|J(t)\|_{L^1}\leq \|J_0\|_{L^1}$;\\
%(iv) $\|J(t)\|_{L^2}\leq \|J_0\|_{L^2}$. ''
%Check whether these are true in our case. 

Denote $x_0$ by the point such that 
\[B_0(x_0)=\max_{x\in\mathbb R} B_0(x).\]
Thanks to the assumptions on $B_0$, we note $0<x_0<1$ and
\[B'_0(x)\geq 0 \ \ \mbox{for} \ x\in(-x_0,x_0);  \ \ B'_0(x)\leq 0 \ \ \mbox{for} \ x\in(-\infty, -x_0)\cap(x_0, \infty).\]
We also note $-x_0$ is the minimum point for $B_0(x)$. 
Let $X(x_0, t)$ be the trajectory of the maximum point of the solution $B(x,t)$ to (\ref{e1d3}) with $\mu=0$ and the initial data $B_0$. It follows from the properties of transport equation that
\begin{equation}\label{traj1}
\frac{d}{dt} X(x_0, t)=J(X(x_0, t), t)=-HB_x(X(x_0, t), t)=-\Lambda B(X(x_0, t), t).
\end{equation}
We note $\Lambda B(X(x_0, t), t)\geq0$ since $X(x_0, t)$ is a maximum point. Hence (\ref{traj1}) implies $X(x_0, t)$ does not increases with respect to time and 
\[0\leq X(x_0, t)\leq x_0<1.\]
Moreover we know
$X(x_0, t)>0$ for $t\in[0,T)$ and
\begin{equation}\label{traj2}
\begin{cases}
B_x(x,t)\geq 0 \ \ \mbox{for} \ x\in(-X(x_0, t),X(x_0, t));  \\
B_x(x,t)\leq 0 \ \ \mbox{for} \ x\in(-\infty, -X(x_0, t))\cap(X(x_0, t), \infty).
\end{cases}
\end{equation}
By the definition of the Hilbert transform we have from (\ref{traj1}) that
\begin{equation}\label{traj3}
\begin{split}
\frac{d}{dt} X(x_0, t)=&-\frac{1}{\pi}P.V.\int_{-\infty}^\infty\frac{B_y(y,t)}{X(x_0, t)-y}\, dy\\
=&-\frac{1}{\pi} \int_{-1}^{-X(x_0, t)}\frac{B_y(y,t)}{X(x_0, t)-y}\, dy-\frac{1}{\pi} P.V.\int_{-X(x_0, t)}^{X(x_0, t)}\frac{B_y(y,t)}{X(x_0, t)-y}\, dy\\
&-\frac{1}{\pi} P.V.\int_{X(x_0, t)}^1\frac{B_y(y,t)}{X(x_0, t)-y}\, dy\\
=&: K_1+K_2+K_3.
\end{split}
\end{equation}
For $-1<y<-X(x_0,t)$, we have $X(x_0,t)-y>2X(x_0,t)>0$. Invoking (\ref{traj2}) we deduce
\begin{equation}\notag
\begin{split}
K_1<&-\frac{1}{2\pi X(x_0,t)} \int_{-1}^{-X(x_0, t)} B_y(y,t)\, dy\\
=&-\frac{1}{2\pi X(x_0,t)}\left(B(-X(x_0,t),t)-B(-1,t)\right)\\
=&\ \frac{1}{2\pi X(x_0,t)}.
\end{split}
\end{equation}
On the other hand, if $-X(x_0,t)\leq y\leq X(x_0,t)$, we have
\[0\leq X(x_0,t)-y\leq 2X(x_0,t).\]
Applying (\ref{traj2}) again gives
\begin{equation}\notag
\begin{split}
K_2\leq& -\frac{1}{2\pi X(x_0,t)} \int_{-X(x_0, t)}^{X(x_0, t)}B_y(y,t)\, dy\\
=&-\frac{1}{2\pi X(x_0,t)}\left(B(X(x_0,t), t)-B(-X(x_0,t), t)\right)\\
=&-\frac{1}{\pi X(x_0,t)}.
\end{split}
\end{equation}
It follows from (\ref{traj2}) that $K_3\leq 0$. Combining the estimates of $K_1$ and $K_2$ with (\ref{traj3}) we get
\begin{equation}\notag
\frac{d}{dt} X(x_0, t)\leq -\frac{1}{2\pi X(x_0,t)}
\end{equation}
which gives
\begin{equation}\notag
\frac{d}{dt} \frac{1}{X(x_0, t)}\geq \frac{1}{\pi X^3(x_0,t)}.
\end{equation}
Hence $\frac{1}{X(x_0, t)}$ becomes unbounded after a short time. It implies $\|B_x(0,t)\|_{L^\infty}$ blows up after a short time. This closes the proof of Theorem \ref{thm-s1}.

\cbdu

\bigskip

%\section{Pure transport case}
%\label{sec-transport}

\section*{Acknowledgement}
The author is indebted to Diego C\'ordoba and Hongjie Dong for their valuable suggestions which have improved significantly the early version of the article.

%\Endrefs
\end{document}